\documentclass[a4paper,11pt,reqno]{smfart}
\usepackage{amssymb,amsmath,mathrsfs,graphics,graphicx,mathptm,float}
\usepackage[T1]{fontenc}
\usepackage[english, francais]{babel}

\theoremstyle{plain}
\newtheorem{thm}{Th\'eor\`eme}[section]
\newtheorem{pro}[thm]{Proposition}
\newtheorem{lem}[thm]{Lemme}
\newtheorem{cor}[thm]{Corollaire}
\newtheorem{theoreme}{Th\'eor\`eme}
\newtheorem{proposition}[theoreme]{Proposition}

\theoremstyle{definition}
\newtheorem*{defi}{D\'efinition}

\newtheorem*{eg}{Exemple}
\newtheorem*{egs}{Exemples}
\newtheorem{rem}[thm]{Remarque}
\newtheorem{rems}[thm]{Remarques}

\def\og{\leavevmode\raise.3ex\hbox{$\scriptscriptstyle\langle\!\langle$~}}
\def\fg{\leavevmode\raise.3ex\hbox{~$\!\scriptscriptstyle\,\rangle\!\rangle$}}

\def\noteA#1#2#3{{\begin{small}$#1$\end{small}} & {\begin{small}#2\end{small}} \hfill{\begin{small}\pageref{#3}\end{small}} \\}

\addtocounter{section}{0}             
\numberwithin{equation}{section}       

\begin{document}
\title{G\'eom\'etrie classique de certains feuilletages quadratiques}

\date{}
\author{D. \textsc{Cerveau}}

\address{Membre de l'Institut Universitaire de France.
IRMAR, UMR 6625 du CNRS, Universit\'e de Rennes $1,$ $35042$ Rennes, France.  
Membre de l'ANR BLAN$06$-$3$\_$137237$}
\email{dominique.cerveau@univ-rennes1.fr}

\author{J. \textsc{D\'eserti}}

\address{Institut de Math\'ematiques de Jussieu, Universit\'e Paris $7,$ Projet G\'eom\'etrie et Dynamique, Site Chevaleret, Case $7012,$ $75205$ Paris Cedex 13, France. 
Membre de l'ANR BLAN$06$-$3$\_$137237$}
\email{deserti@math.jussieu.fr}

\author{D. \textsc{Garba Belko}}
\address{Facult\'e des Sciences, Universit\'e Abdou Moumouni, B.P. 10662 Niamey, Niger.}
\email{garbabelkodjibrilla@yahoo.fr}

\author{R. \textsc{Meziani}}

\address{D\'epartement de Math\'ematiques, Facult\'e des Sciences, Universit\'e Ibn Tofail, Kenitra, Maroc.}
\email{rmeziani@yahoo.com}

\maketitle{}

\begin{altabstract}
\selectlanguage{english}
The set $\mathscr{F}(2;2)$ of quadratic foliations on the complex
projective plane can be identified with a \textsc{Zariski}'s open
set of a projective space of dimension $14$ on which acts 
 $\mathrm{Aut}(\mathbb{P}^2(\mathbb{C})).$ We classify, up to
automorphisms of $\mathbb{P}^2(\mathbb{C}),$ quadratic 
foliations with only one singularity. There are only four
such foliations up to conjugacy; whereas three of them have a dynamic
which can be easily described the dynamic of the fourth is 
still mysterious. This classification also allows us to describe the
action of $\mathrm{Aut}(\mathbb{P}^2(\mathbb{C}))$ on 
$\mathscr{F}(2;2).$ On the one hand we show that the dimension
of the orbits is more than $6$ and that there are exactly two
orbits of dimension $6;$ on the other hand we obtain that
the closure of the generic orbit in $\mathscr{F}
(2;2)$ contains at least seven orbits of dimension~$7$ and
exactly one orbit of dimension $6.$ 

\end{altabstract}

\selectlanguage{french}

\begin{abstract}
L'ensemble $\mathscr{F}(2;2)$ des feuilletages quadratiques
du plan projectif complexe s'identifie \`a un ouvert de \textsc{Zariski} dans
un espace projectif de dimension $14$ sur lequel agit le groupe~$\mathrm{Aut}(\mathbb{P}^2(\mathbb{C})).$
Nous classifions, \`a automorphisme de $\mathbb{P}^2(\mathbb{C})$ pr\`es, 
les feuilletages quadratiques ayant une 
unique singularit\'e. \`A automorphisme pr\`es il y a $4$ feuilletages
ayant cette propri\'et\'e; alors que trois ont une dynamique que l'on
peut d\'ecrire facilement, celle du quatri\`eme reste myst\'erieuse.
Cette classification intervient dans la description de l'action de 
$\mathrm{Aut}(\mathbb{P}^2(\mathbb{C}))$ sur $\mathscr{F}(2;2).$ 
Nous montrons d'une part que la dimension des orbites est 
sup\'erieure ou \'egale \`a $6$ et qu'il y a exactement deux orbites  de
dimension $6$ dont l'une correspond \`a un feuilletage ne pr\'esentant 
qu'un seul point singulier; d'autre part nous obtenons que
l'adh\'erence de l'orbite d'un \'el\'ement g\'en\'erique de $\mathscr{F}(2;2)$ 
contient au moins sept orbites de 
dimension $7$ et une seule orbite de dimension $6.$
\noindent{\it Classification math\'ematique par sujets (2000). --- 37F75, 32S65, 32M25, 32M05, 14L35.}
\end{abstract}


\section*{Introduction}

\noindent Soit $\mathcal{F}$ un feuilletage holomorphe de codimension $1$ et de
degr\'e $N$ sur $\mathbb{P}^2(\mathbb{C}).$ Si $\pi\colon\mathbb{C}^3\setminus\{0\}\to\mathbb{P}^2(\mathbb{C})$ est la projection canonique, le feuilletage homog\`ene $\pi^{-1}\mathcal{F}$ s'\'etend \`a $\mathbb{C}^3$ et est d\'efini par une $1$-forme
$$\omega=a(X,Y,Z)\mathrm{d}X+b(X,Y,Z)\mathrm{d}Y+c(X,Y,Z)\mathrm{d}Z,$$
o\`{u} $a,$ $b$ et $c$ sont des polyn\^{o}mes homog\`enes de
degr\'e $N + 1$ sans composante commune v\'erifiant
l'identit\'e d'\textsc{Euler}: $aX+bY+cZ=0;$ c'est le th\'eor\`eme de \textsc{Chow} pour les feuilletages. Le lieu singulier
$\mathrm{Sing}(\mathcal{F})$\label{not0} de~$\mathcal{F}$ est donn\'e par
$$\pi (\{a=b=c=0\}\setminus \{0\}).$$
\noindent L'identit\'e d'\textsc{Euler} assure 
l'existence de polyn\^omes homog\`{e}nes $p,$ $q,$ $r$ de degr\'e $N$ tels que
$$\omega = p(X,Y,Z) (Y\mathrm{d}X-X\mathrm{d}Y) + 
q(X,Y,Z) (Y\mathrm{d}Z-Z\mathrm{d}Y)+ r(X,Y,Z)(X\mathrm{d}Z-Z\mathrm{d}X).$$
Dans la carte affine $Z=1$ nous dirons abusivement que la $1$-forme $\omega$ s'\'ecrit
$$p(x,y,1)(y\mathrm{d}x-x\mathrm{d}y)-q(x,y,1)\mathrm{d}y-r(x,y,1)\mathrm{d}x.$$

\noindent Rappelons la notion de nombre
de \textsc{Milnor}\label{not1} $\mu(\mathcal{F},m)$ d'un feuilletage
$\mathcal{F}$ en un point singulier~$m.$ Fixons
une carte locale $(u,v)$ telle que $m=(0,0);$
le germe de $\mathcal{F}$ en $m$ est d\'efini, \`a multiplication
par une unit\'e en $0$ pr\`es, par une $1$-forme $E\mathrm{d}u+ F
\mathrm{d}v.$ D\'esignons par $\langle E,F\rangle$\label{not1b} l'id\'eal engendr\'e
par $E$ et $F,$ alors
$$\mu(\mathcal{F},m) = \dim\frac{\mathbb{C}\{u,v\}}{\langle E,F\rangle};$$
c'est aussi la multiplicit\'e d'intersection
$(\mathcal{C}. \mathcal{C}')_0$ des germes de courbes $\mathcal{C}
=(E=0)$ et~$\mathcal{C}'=(F=0).$

\noindent Sur $\mathbb{P}^2(\mathbb{C})$ on dispose d'un th\'eor\`eme de type \textsc{Bezout} (\emph{voir}
\cite{S}); si $\mathcal{F}$ est un feuilletage de degr\'e $N$ on a
\begin{equation}\label{bezout}
\sum_{m \in\mathrm{Sing}(\mathcal{F})} \mu(\mathcal{F}, m)=N^2+N+1.
\end{equation}
Cette formule implique, en particulier, qu'il n'y a pas de feuilletage
r\'egulier ({\it i.e.} sans singularit\'e) sur $\mathbb{P}^2(\mathbb{C}).$  

\noindent Nous noterons 
$\mathscr{F}(2;N)$\label{not7} l'ensemble des feuilletages de degr\'e~$N$ sur $\mathbb{P}^2(\mathbb{C}).$

\noindent La classification des feuilletages de degr\'e $0$ ou $1$ sur
le plan projectif complexe est connue depuis le~XIX$^{\text{\`eme}}$ si\`ecle~(\cite{J}): un feuilletage de degr\'e
$0$ sur $\mathbb{P}^2(\mathbb{C})$ est un pinceau de droites et poss\`ede une seule singularit\'e; tout
feuilletage de degr\'e $1$ sur le plan projectif complexe poss\`ede trois
singularit\'es compt\'ees avec multiplicit\'e, a, au moins, une droite
invariante et est donn\'e par une forme ferm\'ee rationnelle. Pour~$N\geq2$ peu de 
propri\'et\'es ont \'et\'e \'etablies, hormis la non existence g\'en\'erique de courbe
invariante~(\cite{J, CLN}). En particulier le probl\`eme du \og minimal 
exceptionnel\fg\hspace{1mm} n'est pas r\'esolu, m\^eme en degr\'e~$2$ (\emph{voir} 
\cite{CLNS}). Ce probl\`eme, qui va dans le sens d'un \'enonc\'e de type Poincar\'e-Bendixon, consiste
\`a d\'ecrire les diff\'erentes possibilit\'es pour les adh\'erences (ordinaires) des feuilles. En particulier on 
ignore actuellement si une telle adh\'erence doit contenir n\'ecessairement un point singulier. Une 
feuille dont l'adh\'erence ne contiendrait pas de singularit\'e produirait alors un \og minimal
exceptionel\fg. 
Il y a a priori plusieurs fa\c{c}ons
de l'aborder; soit on cherche \`a d\'egager des
propri\'et\'es des feuilletages g\'en\'eriques, soit au contraire on \'etudie les
d\'eg\'en\'erescences les plus compliqu\'ees que l'on essaie
ensuite de d\'eformer. C'est dans cette derni\`ere optique que
nous nous proposons de classifier, \`a isomorphisme pr\`es, les
feuilletages quadratiques de $\mathbb{P}^2(\mathbb{C})$ ayant une unique
singularit\'e que l'on suppose \^etre le point~$0$ de coordonn\'ees
$(0,0)$ dans la carte affine $Z=1.$ On peut a priori penser que ces feuilletages, du fait de la confluence de leurs singularit\'es en une seule, sont tr\`es compliqu\'es. Comme nous le verrons cette philosophie est dans certains cas mise en d\'efaut. La formule~(\ref{bezout})
implique que $\mu(\mathcal{F},0)=7.$ En particulier les germes de $P$
et $Q$ en $0$ ne forment pas un syst\`eme de coordonn\'ees local en $0;$ en
effet si c'\'etait le cas $\mu(\mathcal{F},0)$ vaudrait $1.$  Par suite, \`a isomorphisme lin\'eaire 
pr\`es, nous avons les trois possibilit\'es suivantes
\begin{itemize}
\item[\texttt{1. }] $\omega=y\mathrm{d}y+
\text{termes de plus haut degr\'e},$ la singularit\'e est dite de type {\sl nilpotent}, \'eventualit\'e \'etudi\'ee au \S\hspace{1mm}\ref{sectionnilpotent};

\item[\texttt{2. }] $\omega=x\mathrm{d}y+
\text{termes de plus haut degr\'e},$ la singularit\'e $0$ est dite de type {\sl selle-noeud}, 
ce cas sera l'objet du \S\hspace{1mm}\ref{sectionsellenoeud};

\item[\texttt{3. }] le $1$-jet de $\omega$ est nul en $0$ (\emph{voir} \S\hspace{1mm}\ref{section1jetnul}).
\end{itemize}

\bigskip

\noindent En analysant le point singulier $0,$ plus pr\'ecis\'ement 
en traduisant dans chacune des \'eventualit\'es qui pr\'ec\`edent,
l'\'egalit\'e $\mu(\mathcal{F},0)=7$ nous obtenons la description
des feuilletages quadratiques de $\mathbb{P}^2(\mathbb{C})$
ayant une unique singularit\'e.

\begin{theoreme}\label{feuilquad}
{\sl \`A automorphisme de $\mathbb{P}^2(\mathbb{C})$ pr\`es, il y a
quatre feuilletages quadratiques $\mathcal{F}_1,$ $\ldots,$ $\mathcal{F}_4$ sur le plan projectif complexe ayant une
seule singularit\'e. Ils sont d\'ecrits respectivement en carte affine par les $1$-formes
suivantes

\begin{itemize}
\item[\texttt{1. }] $\omega_1=x^2\mathrm{d}x+y^2(x\mathrm{d}y-y\mathrm{d}x);$\label{not2}

\item[\texttt{2. }] $\omega_2=x^2\mathrm{d}x+(x+y^2)(x\mathrm{d}y-y\mathrm{d}x);$\label{not3}

\item[\texttt{3. }] $\omega_3=xy\mathrm{d}x+(x^2+y^2)(x\mathrm{d}y-y\mathrm{d}x);$\label{not4}

\item[\texttt{4. }] $\omega_4=(x+y^2-x^2y)\mathrm{d}y+x(x+y^2)\mathrm{d}x.$\label{not5}
\end{itemize}}
\end{theoreme}

\noindent On observe qu'aucun de ces feuilletages pr\'esente de singularit\'e
nilpotente. Notons que le feuilletage d\'efini par les niveaux de 
$\frac{x+y^2}{x^2}$ est de degr\'e $1$ et pr\'esente une unique singularit\'e qui de plus est nilpotente.

\noindent 
Les 
feuilletages~$\mathcal{F}_1,$ $\mathcal{F}_2$ et $\mathcal{F}_3$ admettent
respectivement pour int\'egrale premi\`ere 
\begin{small}
\begin{align*}
& \frac{1}{3}\left(\frac{y}{x}\right)^3-\frac{1}{x}, &&\left(2+\frac{1}{x}+
2\left(\frac{y}{x}\right)+\left(\frac{y}{x}\right)^2\right)\exp\left(-\frac{y}{x}\right) &&
\text{et}&&\left(\frac{y}{x}\right)\exp\left(\frac{1}{2}\left(\frac{y}{x}\right)^2-\frac{1}{x}\right).
\end{align*}
\end{small}
\noindent Notons que l'adh\'erence d'une feuille g\'en\'erique de
$\mathcal{F}_1$ est une cubique cuspidale; remarquons aussi que~$\mathcal{F}_1,$ $\mathcal{F}_2$ et $\mathcal{F}_3$ peuvent \^etre d\'efinis par une $1$-forme
ferm\'ee rationnelle. Chacun de ces trois mod\`eles compte au moins une
courbe alg\'ebrique invariante; il n'en est pas de m\^eme pour $\mathcal{F}_4.$ En
particulier $\mathcal{F}_4$ ne peut poss\'eder de structure
transversalement projective~(\cite{CLNLPT}); ceci implique d'ailleurs
que $\mathcal{F}_4$ n'a pas d'int\'egrale premi\`ere de type \textsc{Liouville}.
Les int\'egrales premi\`eres de
$\mathcal{F}_1,$ $\mathcal{F}_2,$ $\mathcal{F}_3$ nous 
permettent de d\'ecrire leurs feuilles de fa\c{c}on imm\'ediate.

\begin{proposition}
{\sl Les feuilletages $\mathcal{F}_1,$ $\mathcal{F}_2$ et $\mathcal{F}_3$
poss\`edent la propri\'et\'e suivante: pour tout point~$m$ r\'egulier
la feuille de $\mathcal{F}_i$ passant par $m$ est une courbe 
enti\`ere, {\it i.e.} une courbe param\'etr\'ee par une application holomorphe de 
$\mathbb{C}$ dans $\mathbb{C}^2.$

\noindent En fait $\mathcal{F}_1,$ $\mathcal{F}_2$ et $\mathcal{F}_3$
sont d\'efinis, dans la carte $X=1,$ par des champs de vecteurs polynomiaux
complets.}
\end{proposition}

\noindent Nous n'avons aucune
id\'ee de la nature des feuilles de $\mathcal{F}_4;$
par exemple on peut se demander si $\mathcal{F}_4$
a un minimal exceptionnel (\cite{CLNS}), {\it i.e.} s'il existe 
une feuille de $\mathcal{F}_4$ n'adh\'erant pas
\`a la singularit\'e. L'\'etude de ce probl\`eme sur cet exemple
nous semble pertinent. Dans \cite{CdF} \textsc{Camacho} et \textsc{de Figueiredo}
montrent par des techniques num\'eriques que le feuilletage de \textsc{Jouanolou}
de degr\'e $2$ (nous en reparlerons plus loin) ne poss\`ede pas de minimal exceptionnel.
Il serait int\'eressant d'adapter leurs m\'ethodes pour l'\'etude du feuilletage $\mathcal{F}_4.$

\bigskip

\noindent Un th\'eor\`eme de \textsc{Luna} et \textsc{Vust} 
affirme que si $\mathrm{M}$ est une vari\'et\'e alg\'ebrique affine, $\mathrm{G}$ un 
groupe r\'eductif agissant alg\'ebriquement sur $\mathrm{M}$ avec isotropie
g\'en\'erique r\'eductive, l'orbite g\'en\'erique de $\mathrm{G}$
est ferm\'ee dans $\mathrm{M}$ (\emph{voir} \cite{LV}). Dans cet esprit nous 
nous demandons si l'orbite $\mathcal{O}(\mathcal{F})$ sous l'action 
du grou\-pe~$\mathrm{Aut}(\mathbb{P}^2(\mathbb{C}))=\mathrm{PGL}_3
(\mathbb{C})$ d'un \'el\'ement g\'en\'erique $\mathcal{F}$ de 
$\mathscr{F}(2;2)$ est ferm\'ee dans $\mathscr{F}(2;2).$ La r\'eponse 
est non et la classification pr\'ec\'edente entre ici en jeu.

\begin{theoreme}\label{degenerea}
{\sl Il existe un ensemble alg\'ebrique $\Sigma$ non trivial, contenu dans 
$\mathscr{F}(2;2)$ ayant la propri\'et\'e suivante: pour tout feuilletage 
$\mathcal{F}$ de $\mathscr{F}(2;2)\setminus\Sigma$ l'adh\'erence de l'orbite 
$\mathcal{O}(\mathcal{F})$ de $\mathcal{F}$ contient $\mathcal{F}_1.$ 

\noindent En particulier, pour tout $\mathcal{F}$ dans $\mathscr{F}(2;2)
\setminus\Sigma,$ l'orbite $\mathcal{O}(\mathcal{F})$ de $\mathcal{F}$ n'est pas
ferm\'ee.}
\end{theoreme}

\noindent Donnons une interpr\'etation g\'eom\'etrique de cet ensemble $\Sigma.$
Un point r\'egulier $m$ du feuilletage $\mathcal{F}$ est dit d'inflexion ordinaire
 si la feuille passant par~$m$ a un point 
d'inflexion ordinaire en $m;$ d\'esignons par $\mathrm{Flex}(\mathcal{F})$
l'adh\'erence de ces points. Un feuilletage~$\mathcal{F}$ appartient \`a 
$\mathscr{F}(2;2)\setminus\Sigma$ si et seulement si $\mathrm{Flex}
(\mathcal{F})$ est non vide. \medskip

\noindent Le Th\'eor\`eme \ref{degenerea} n'est pas en contradiction 
avec l'\'enonc\'e de \textsc{Luna} et \textsc{Vust} puisque 
$\mathscr{F}(2;2)$ n'est pas une vari\'et\'e affine.\medskip

\noindent La dimension de l'orbite de $\mathcal{F}_1$ est $6$ qui est 
la dimension minimale possible, et ce en tout
degr\'e sup\'erieur ou \'egal \`a $2$ (Proposition \ref{dimmin}). Nous 
montrons qu'il y a exactement deux orbites
de dimension $6$ (Proposition~\ref{dim6}); la seconde est associ\'ee au feuilletage~$\mathcal{F}_5$
donn\'e en carte affine par 
$$\omega_5=x^2\mathrm{d}y+y^2(x\mathrm{d}y-y\mathrm{d}x).$$\label{not5b}

\noindent Ce feuilletage poss\`ede deux points singuliers et l'int\'egrale premi\`ere
$\frac{y}{x}-\frac{1}{y};$ autrement dit il s'agit du pinceau de coniques $(y^2-xz)+\lambda 
xy.$ Pour des raisons de dimension les orbites $\mathcal{O}(\mathcal{F}_1)$
et $\mathcal{O}(\mathcal{F}_5)$ sont ferm\'ees.

\noindent De nombreux probl\`emes, en particulier la partie b) du
$16^{\text{\`eme}}$ probl\`eme de \textsc{Hilbert}, sont abord\'es par 
des techniques de perturbation en g\'en\'eral \`a partir de feuilletages
\og hamiltoniens\fg. Ce qui pr\'ec\`ede montre en un certain sens
que presque tout feuilletage quadratique est une \og petite 
perturbation\fg\hspace{1mm} du feuilletage \og hamiltonien\fg\hspace{1mm}
$\mathcal{F}_1.$ Nous pr\'eciserons le Th\'eor\`eme \ref{degenerea}, en 
particulier dans l'\'enonc\'e qui suit.
 
\begin{theoreme}
{\sl Si $\mathcal{F}$ d\'esigne un \'el\'ement g\'en\'erique de $\mathscr{F}(2;2)$
l'adh\'erence de l'orbite $\mathcal{O}(\mathcal{F})$ de~$\mathcal{F}$ contient 
au moins sept orbites de dimension~$7$ et une seule orbite de dimension~$6,$
celle de~$\mathcal{F}_1.$}
\end{theoreme}

%
%
%

\noindent D'autres consid\'erations sur les adh\'erences d'orbites seront 
abord\'ees dans le texte.

\bigskip

\noindent\textit{\textbf{Remerciements.}} Merci \`a L. \textsc{Pirio} 
pour sa disponibilit\'e et \`a J. \textsc{Pereira} pour ses remarques.

\noindent Les deux premiers auteurs remercient l'IUF et
l'ANR Symplexe (ANR BLAN$06$-$3$\_$137237$) qui ont contribu\'e au 
bon d\'eroulement de cette collaboration ainsi que le CIRM pour les excellentes
conditions de travail dont ils ont b\'en\'efici\'e.

\noindent Les troisi\`eme et quatri\`eme auteurs remercient l'IRMAR pour leurs s\'ejours
\`a Rennes.

\section{Feuilletages quadratiques de $\mathbb{P}^2(
\mathbb{C})$ ayant une seule singularit\'e}

\noindent Consid\'erons un feuilletage quadratique $\mathcal{F}$ sur $\mathbb{P}^2
(\mathbb{C})$ ayant une unique singularit\'e;  la classification menant au 
Th\'eor\`eme \ref{feuilquad} est \'etablie au cas par cas
suivant la nature du~$1$-jet de $\omega$ d\'efinissant $\mathcal{F}$ au point 
singulier que nous supposons \^etre l'origine $0$ de la carte
affine $Z=1.$ On \'ecrira $\omega$ sous la forme $A\mathrm{d}x+B\mathrm{d}y+\phi
(x\mathrm{d}y-y\mathrm{d}x)$ et on adoptera les notations suivantes
\begin{align*}
&A=a_{1,0}x+a_{0,1}y+a_{2,0}x^2+a_{1,1}xy+a_{0,2}y^2,&&
B=b_{1,0}x+b_{0,1}y+b_{2,0}x^2+b_{1,1}xy+b_{0,2}y^2,
\end{align*}
$$\phi=\phi_{2,0}x^2+\phi_{1,1}xy+\phi_{0,2}y^2.$$ On d\'esignera par $P$ (resp. $Q$) la composante en $\mathrm{d}x$
(resp. $\mathrm{d}y$); autrement dit on pose $P=A-y\phi$ et $Q=B+x\phi.$

\subsection{\'Etude du cas nilpotent}\label{sectionnilpotent}

\begin{pro}\label{casnil}
{\sl Soit $\mathcal{F}$ un feuilletage quadratique sur
$\mathbb{P}^2(\mathbb{C})$ dont le lieu singulier est r\'eduit \`a un 
point; cette singularit\'e ne peut \^etre nilpotente.}
\end{pro}

\begin{proof}[\sl D\'emonstration]
Raisonnons par l'absurde: supposons que $\mathcal{F}$ poss\`ede
une unique singularit\'e et que celle-ci soit nilpotente. \`A
isomorphisme lin\'eaire pr\`es $\mathcal{F}$ est d\'efini par  la
$1$-forme~$\omega=P\mathrm{d}x+Q\mathrm{d}y$ o\`u 
\begin{align*}
&P(x,y)=a_{2,0} x^2+a_{0,2}y^2+a_{1,1}xy-y(\phi_{2,0}x^2+\phi_{0,2}y^2+\phi_{1,1}xy),\\
&Q(x,y)=y+b_{2,0} x^2+b_{0,2}y^2+b_{1,1}  xy+x(\phi_{2,0}x^2+\phi_{0,2}y^2+\phi_{1,1}xy).
\end{align*}

\noindent Remarquons que la condition $\mu(\mathcal{F},0)=7$ implique la
nullit\'e de $a_{2,0}.$ Dans le cas contraire $\langle P,Q\rangle_{,0}=\langle y,x^2\rangle_{,0}$
dans l'anneau local $\mathbb{C}\{x,y\}$ qui implique $\mu(\mathcal{F},0)
=2.$\medskip

\noindent Ainsi le polyn\^ome $P$ est divisible par~$y$
\begin{align*}
& P=yP_1 && \text{avec} && P_1=a_{0,2}y+	a_{1,1}
x-(\phi_{2,0}x^2+\phi_{0,2}y^2+\phi_{1,1}xy).
\end{align*}
Comme $\omega$ est \`a singularit\'e isol\'ee en $0,$ le polyn\^ome $Q$ est non divisible
par $y$ et 
$\vert b_{2,0} \vert + \vert \phi_{2,0} \vert$ est non nul.

\noindent Posons $\mathcal{C}=(Q=0)$ et $\mathcal{C'}=(P=0);$ nous avons
$$\mu (\mathcal{F}, 0)= (\mathcal{C}\cdot\mathcal{C}')_0 = (\mathcal{C}\cdot(y=0))_0 + (\mathcal{C}\cdot(P_1
=0))_0.$$
\noindent Notons que la multiplicit\'e d'intersection
$(\mathcal{C}\cdot(y=0))_0$ vaut $2$ ou $3.$

\noindent Si $a_{1,1}$
est non nul~$(P_1 = 0)$ est transverse \`a $\mathcal{C}$ et
$(\mathcal{C}\cdot(P_1=0))_0 = 1;$ d'o\`u l'in\'egalit\'e $\mu
(\mathcal{F}, 0)\leq 4:$ contradiction. D'o\`u l'\'egalit\'e $a_{1,1}=0.$\medskip

\noindent Sur la droite $y=0$ la $1$-forme $\omega$ s'\'ecrit
$x^2(b_{2,0}+\phi_{2,0}x)\mathrm{d}y;$ si $\phi_{2,0}b_{2,0}\not=0$ le point
$(-\frac{b_{2,0}}{\phi_{2,0}},0)$ est une singularit\'e de $\omega,$
distincte de $(0,0),$ ce qui est, par hypoth\`ese, impossible. Il
en r\'esulte que le produit $\phi_{2,0}b_{2,0}$ est nul. \medskip

\noindent Montrons, par l'absurde, que $\phi_{2,0}\neq 0.$ Si $\phi_{2,0}=0,$ alors d'apr\`es ce 
qui pr\'ec\`ede~$b_{2,0}\not=~0$~et
$$\omega= \left(y+b_{2,0}x^2 +b_{0,2} y^2+b_{1,1}
xy+x(\phi_{0,2}y^2+\phi_{1,1}xy)\right)\mathrm{d}y+y^2\left(a_{0,2} -\phi_{0,2}y-\phi_{1,1}x\right)\mathrm{d}x.$$
\noindent Si $a_{0,2}$ est diff\'erent de $0,$ alors $a_{0,2}
-\phi_{0,2}y-\phi_{1,1}x$ est une unit\'e en $(0,0);$ de sorte que
$$\mu (\mathcal{F}, 0)= 2(\mathcal{C}\cdot(y=0))_0 =4,$$
\noindent ce qui est absurde.
Donc $a_{0,2} =0$ et $\omega$ est du type
$$\left(y+b_{2,0} x^2 + b_{0,2} y^2+b_{1,1} xy+x(\phi_{0,2}y^2+\phi_{1,1}xy)\right)\mathrm{d}y-
y^2\left(\phi_{0,2}y+\phi_{1,1}x\right)\mathrm{d}x.$$
\noindent Si $\phi_{1,1}\neq 0$ nous constatons que $\mu (\mathcal{F}, 0)=5;$ enfin si
$\phi_{1,1}=0$ le coefficient $\phi_{0,2}$ est non nul et~$\mu(\mathcal{F}, 0)=~6,$ ces deux
cas sont exclus. Donc $\phi_{2,0}\neq 0.$

\medskip
\noindent Finalement $\omega$ s'\'ecrit
\begin{small}
\begin{align*}
&\left(y+b_{0,2} y^2+b_{1,1} xy+x(\phi_{2,0}x^2+\phi_{0,2}y^2+\phi_{1,1}xy)\right)\mathrm{d}y+y\left(
a_{0,2} y-(\phi_{2,0}x^2+\phi_{0,2}y^2+\phi_{1,1}xy)\right)\mathrm{d}x, && \phi_{2,0}\neq 0.
\end{align*}
\end{small}
\noindent La multiplicit\'e d'intersection en $0$ entre $(Q=0)$ et
$(y=0)$ est $3.$ En outre nous avons au niveau des id\'eaux dans
l'anneau des germes de fonctions holomorphes en $0$
\begin{eqnarray}
&&\langle Q,a_{0,2} y-(\phi_{2,0}x^2+\phi_{0,2}y^2+\phi_{1,1}xy)\rangle_{,0}\nonumber\\
&&\hspace{6mm}=
\langle y(1+b_{0,2} y+(b_{1,1} +  a_{0,2})x) , a_{0,2}
y-(\phi_{2,0}x^2+\phi_{0,2}y^2+\phi_{1,1}xy)\rangle_{,0}\nonumber\\
&&\hspace{6mm}=\langle y,x^2\rangle_{,0}.\nonumber
\end{eqnarray}

\noindent Nous en d\'eduisons que $\mu (\mathcal{F},0) =3+2=5.$ Il n'existe donc
pas de feuilletage quadratique de~$\mathbb{P}^2(\mathbb{C})$ ayant une
seule singularit\'e qui soit de type nilpotent.
\end{proof}

\subsection{Singularit\'es de type selle-noeud}\label{sectionsellenoeud}

\begin{pro}\label{casselnoeud}
{\sl Soit $\mathcal{F}$ un feuilletage quadratique ayant une seule 
singularit\'e. Si cette singularit\'e est de type selle-n\oe ud, $\mathcal{F}$
est conjugu\'e au feuilletage $\mathcal{F}_4$ d\'ecrit par la $1$-forme
$$\omega_4=(x+y^2-x^2y)\mathrm{d}y+x(x+y^2)\mathrm{d}x.$$}
\end{pro}

\begin{proof}[\sl D\'emonstration]
Supposons comme pr\'ec\'edemment que $0$ soit l'unique
singularit\'e de~$\mathcal{F}.$ Une fois encore nous allons traduire
l'\'egalit\'e $\mu(\mathcal{F},0)=7.$ \`A isomorphisme pr\`es~$\omega$ 
s'\'ecrit $P\mathrm{d}x+Q\mathrm{d}y$ dans la carte~$Z=1$ avec
$$P(x,y)=a_{2,0} x^2 + a_{0,2} y^2 +a_{1,1} xy-y(\phi_{2,0}x^2+\phi_{0,2}y^2
+\phi_{1,1}xy),$$
$$Q(x,y) = x+b_{2,0} x^2+b_{0,2} y^2+b_{1,1}xy +x(\phi_{2,0}x^2+\phi_{0,2}y^2
+\phi_{1,1}xy).$$
Montrons par l'absurde que la droite $(x=0)$ n'est pas invariante par
$\mathcal{F}.$ Si $(x=0)$ est invariante par~$\mathcal{F},$ le coefficient 
$b_{0,2}$ est nul. 
Comme $\omega$ est \`a singularit\'e isol\'ee
$\vert a_{0,2} \vert + \vert \phi_{0,2} \vert$ est non
nulle; on en d\'eduit l'\'egalit\'e d'id\'eaux
$$\langle P,Q\rangle_{,0}=\langle a_{0,2}y^2 +\phi_{0,2}y^3, x\rangle_{,0}$$
qui conduit \`a $\mu(\mathcal{F}, 0) \leq 3:$ contradiction.

\noindent  Consid\'erons un point g\'en\'erique $m_0=(0,y_0)$ 
de l'axe $(x=0)$ et $\mathcal{D}$ la tangente \`a la feuille de~$\mathcal{F}$ 
passant par~$m_0$; cette tangente est donc distincte de
l'axe $(x=0).$ Quitte \`a faire agir un automorphisme de $\mathbb{P}^2(\mathbb{C})$
pr\'eservant les axes $(x=0)$ et~$(y=0),$ envoyant $m_0$ sur le
point \`a l'infini $(0:1:0)$ et~$\mathcal{D}$ sur la droite \`a l'infini
$Z=0,$ nous pouvons nous ramener \`a $\phi_{0,2}=0.$ Maintenant, quitte \`a faire agir la transformation lin\'eaire diagonale $(b_{0,2}x,y),$
on peut effectuer la normalisation $b_{0,2}=1.$

\noindent Le th\'eor\`eme des fonctions implicites assure que la courbe
$Q=0$ est en $0$ un graphe local (calcul r\'ealis\'e avec Maple)
\begin{align*}
& x=x(y) =-y^2+b_{1,1}y^3-(b_{1,1}^2+b_{2,0})y^4+(3b_{2,0}b_{1,1}-\phi_{1,1}
+b_{1,1}^3)y^5\\
&\hspace{1cm}+(\phi_{2,0}-b_{1,1}^4+3\phi_{1,1}b_{1,1}-2b_{2,0}^2-
6b_{2,0}b_{1,1}^2)y^6\\
&\hspace{1cm}+(10b_{2,0}^2b_{1,1}-4b_{2,0}\phi_{1,1}+10b_{2,0}
b_{1,1}^3-4\phi_{2,0}b_{1,1}-6\phi_{1,1}b_{1,1}^2+b_{1,1}^5)y^7\text{ mod } y^8.
\end{align*}

\noindent La condition $\mu(\mathcal{F},0)=7$ se traduit par le fait que 
$P(x(y),y)$ est d'ordre pr\'ecis\'ement $7.$ En \'ecrivant explicitement cette
condition on obtient les \'egalit\'es suivantes
\begin{align*}
&a_{0,2}=a_{1,1}=0, &&\phi_{1,1}=-a_{2,0},&&
\phi_{2,0}=-a_{2,0}b_{1,1}, &&a_{2,0}b_{2,0}=0 && a_{2,0}(3b_{1,1}b_{2,0}
+a_{2,0})\not=0,
\end{align*}
\noindent ces deux derni\`eres \'egalit\'es impliquant~$b_{2,0}=0$ et $a_{2,0}\not=0.$

\noindent Ainsi 
$$\omega = \left(x+y^2+b_{1,1} xy -x^2(b_{1,1} x+y)\right)\mathrm{d}y
+x\left( x+y(b_{1,1} x+y)\right)\mathrm{d}x.$$
Quitte \`a remplacer $\omega$ par $\varphi^\ast\omega$, o\`u
$\varphi(x,y)=\left(\frac{x}{1-b_{1,1} y}, \frac{y}{1-b_{1,1} y}\right),$
le feuilletage $\mathcal{F}$ est d\'efini, dans les coordonn\'ees affines $(x,y),$
par la $1$-forme $$\omega_4=(x+y^2-x^2y)\mathrm{d}y+x(x+y^2)\mathrm{d}x.$$
\end{proof}

\medskip
\noindent Les feuilletages de $\mathbb{P}^2(\mathbb{C})$ les plus
\og simples\fg\, sont ceux qui poss\`edent une structure transversalement 
projective (\cite{CLNLPT}). La proposition qui suit dit en particulier que 
$\mathcal{F}_4$ ne poss\`ede pas de telle structure; il devient alors naturel d'aborder pour ce feuilletage le probl\`eme du minimal exceptionnel. Le portrait de phase \og r\'eel\fg\, de la restriction de $\mathcal{F}_4$ \`a $\mathbb{P}^2(\mathbb{R})$ pr\'esent\'e au \S\,\ref{portraitdephase} fait penser que toutes les feuilles adh\`erent \`a la singularit\'e. 

\begin{pro}
{\sl Le feuilletage $\mathcal{F}_4$ n'admet pas de courbe invariante
alg\'ebrique. En particulier $\mathcal{F}_4$ n'est pas transversalement projectif.}
\end{pro}

\begin{proof}[\sl D\'emonstration]
Rappelons que si $\mathcal{C}$ est une courbe alg\'ebrique invariante par
$\mathcal{F}_4,$ elle
contient au moins une singularit\'e de $\mathcal{F}_4,$ {\it i.e.} 
$\mathcal{C}$ passe n\'ecessairement par l'unique point singulier
$(0:0:1)$ de $\mathcal{F}_4$ (\emph{voir} \cite{CLN, J}). Comme  ce dernier est de type selle-noeud, toute courbe
invariante locale est soit lisse, soit \`a croisement ordinaire, c'est-\`a-dire
du type $xy+\text{termes de degr\'e sup\'erieur}=0$ (\emph{voir} \cite{CM}).

\noindent D'apr\`es \cite{CLN} le degr\'e de $\mathcal{C}$ est inf\'erieur ou \'egal \`a $4.$ Si
$\deg\mathcal{C}=4$ alors $\mathcal{C}$ est r\'eductible et~$\mathcal{F}_4$ est
d\'ecrit par une forme logarithmique $\sum \lambda_i\frac{
d\mathrm{f}_i}{\mathrm{f}_i}$(\emph{voir} \cite{CLN}, th\'eor\`eme $1$); en
particulier le germe de $\mathcal{F}_4$ en $0$ est donn\'e par une forme
logarithmique. Or ce n'est pas le cas pour une singularit\'e de
type selle-noeud (\cite{CM}); donc $\deg\mathcal{C}\leq3.$ Si $\mathcal{C}$ est irr\'eductible lisse
de degr\'e~$3,$ alors (\cite{CLN}, proposition 3) assure que
$\mathcal{F}_4$ a une int\'egrale premi\`ere
rationnelle; il s'en suit que le germe de $\mathcal{F}_4$ en $0$ a une int\'egrale
premi\`ere m\'eromorphe ce qui est encore impossible pour une
singularit\'e selle-noeud. Si $\mathcal{C}$ est irr\'eductible mais non lisse, il s'agit
d'une cubique \`a point double donn\'ee par une \'equation de la forme
$xy+\psi=0,$ o\`u $\psi$ est un polyn\^ome homog\`ene de degr\'e $3.$ Si
$\mathcal{C}$ est r\'eductible elle poss\`ede une branche de degr\'e plus
petit que $2,$ admettant une \'equation de l'un des deux  types
qui suivent
\begin{align*}
& x+ax^2+bxy+cy^2=0
&& \text{ ou }
&&
y+ax^2+bxy+cy^2=0;
\end{align*}
de m\^eme si le degr\'e de $\mathcal{C}$ est inf\'erieur ou \'egal \`a $2.$
Dans toutes ces situations on montre par un calcul Maple que $\omega_4$ ne poss\`ede pas de telles
courbes invariantes.

\noindent D'apr\`es (\cite{CLNLPT}, corollaire 2.16) le feuilletage
$\mathcal{F}_4$ n'est pas transversalement projectif, la pr\'esence d'une telle 
structure n\'ecessitant la pr\'esence d'une courbe alg\'ebrique invariante.
\end{proof}

\begin{rem}
Il r\'esulte des travaux de \textsc{Singer} (\cite{Si}) que
$\mathcal{F}_4$ n'est pas int\'egrable au sens de \textsc{Liouville}. 
De m\^eme $\mathcal{F}_4$ ne peut \^etre d\'efini
par une forme ferm\'ee rationnelle.
\end{rem}

\subsection{Cas 1-jet nul}\label{section1jetnul}\hspace{1mm}

\noindent Dans les trois lemmes qui suivent $\mathcal{F}$ d\'esigne un 
feuilletage quadratique sur $\mathbb{P}^2(\mathbb{C})$ d\'efini par une~$1$-forme $\omega$ et tel que
\begin{itemize}
\item[\texttt{1. }] l'unique singularit\'e de $\mathcal{F}$ soit $(0:0:1);$ 

\item[\texttt{2. }] le $1$-jet en~$0$ de~$\omega$ soit nul. 
\end{itemize}
Dans ce cas
$$\omega=A(x,y)\mathrm{d}x+B(x,y)\mathrm{d}y+\phi(x,y)(x\mathrm{d}y-y\mathrm{d}x),$$
o\`u $A,$ $B$ et~$\phi$ sont des polyn\^omes homog\`enes de degr\'e
$2.$ Le feuilletage $\mathcal{F}$ \'etant quadratique, le c\^one tan\-gent~$xA+yB$ de $\omega$ en $0$ ne peut
\^etre identiquement nul. Le polyn\^ome $\phi$
n'est pas non plus identiquement nul sinon la droite \`a l'infini
serait invariante par~$\mathcal{F}$ qui poss\`ederait donc une
singularit\'e sur cette droite, ce qui est exclu. Nous allons raisonner
suivant la nature du c\^one tangent qui, a priori, peut \^etre trois 
droites, deux droites ou une droite.

\begin{lem}\label{1jetnul1}
{\sl Le c\^one tangent de $\omega$ ne peut \^etre l'union de $3$ droites distinctes.}
\end{lem}

\begin{proof}[\sl D\'emonstration]
Raisonnons par l'absurde; nous pouvons supposer que le c\^one tangent
est form\'e des droites $x=0,$ $y=0$ et $y-x=0;$ il s'en suit que~$\omega$ s'\'ecrit (\cite{CM})
\begin{small}
\begin{align*}
& xy(y-x)\left(\lambda_1\frac{\mathrm{d}x}{x}+\lambda_2\frac{\mathrm{d}y}{y}+
\lambda_3\frac{\mathrm{d}(y-x)}{y-x}\right)
+\left(\phi_{2,0} x^2+\phi_{1,1} xy+\phi_{0,2} y^2\right)(x\mathrm{d}y-y\mathrm{d}x), && \lambda_i,\,\phi_{i,j}\in\mathbb{C}.
\end{align*}
\end{small}
\noindent Sur la droite $x=0$ nous avons $\omega=y^2(\lambda_1-\phi_{0,2} y)\mathrm{d}x.$ 
Comme $0$ est l'unique singularit\'e de $\mathcal{F},$ le pro\-duit~$\lambda_1\phi_{0,2}$ est nul. De la m\^eme fa\c{c}on, en restreignant
$\omega$ aux droites $y=0$ et $y-x=0$ nous obtenons
$$\lambda_2\phi_{2,0}=\lambda_3(\phi_{2,0}+\phi_{1,1}+\phi_{0,2})=0.$$
Le point $(1:0:0)$ \'etant non singulier nous avons $\phi_{2,0}\neq0$  et par sui\-te~$\lambda_2=0.$ Un argument analogue montre que $\lambda_1$ et
$\lambda_3$ sont nuls, mais ceci contredit l'\'egalit\'e
$\deg\mathcal{F}=2.$
\end{proof}

\begin{lem}\label{1jetnul2}
{\sl Si le c\^one tangent de $\omega$ en $0$ est 
compos\'e de deux droites distinctes, alors, \`a isomorphisme pr\`es, 
$\mathcal{F}$ est d\'efini par
$$\omega_3=xy\mathrm{d}x+(x^2+y^2)(x\mathrm{d}y-y\mathrm{d}x).$$}
\end{lem}

\begin{proof}[\sl D\'emonstration]
\`A conjugaison pr\`es $\mathcal{F}$ est donn\'e par la
$1$-forme (\cite{CM})
\begin{align*}
&\omega=x^2y\left(
\lambda\frac{\mathrm{d}x}{x}+\delta\frac{\mathrm{d}y}{y}+\varepsilon \mathrm{d}\left(\frac{y}{x}\right)\right)+\left(\phi_{2,0} 
x^2+\phi_{1,1} xy+\phi_{0,2} y^2\right)(x\mathrm{d}y-y\mathrm{d}x),&&\lambda,\,
\delta,\,\varepsilon,\,\phi_{i,j}\in\mathbb{C}.
\end{align*}
Comme au Lemme \ref{1jetnul1}, en consid\'erant les restrictions de~$\omega$ 
aux droites~$x=0$ et~$y=0,$ nous \'etablissons que 
$$\varepsilon\phi_{0,2}=\delta\phi_{2,0}=0.$$
Le point $(1:0:0)$ n'est pas singulier donc
$\phi_{2,0}\neq0$ et $\delta=0.$ De m\^eme en \'ecrivant 
que $(0:1:0)$ n'est pas singulier nous remarquons
que $\phi_{0,2}\neq0$ et par suite que $\varepsilon=0.$ Le feuilletage $\mathcal{F}$
\'etant quadratique $\lambda$ est non nul; nous pouvons donc supposer
que~$\lambda=1,$ d'o\`u
$$\frac{\omega}{x^3y}=\frac{\mathrm{d}x}{x^2}+\frac{(\phi_{2,0} x^2+\phi_{1,1} xy+\phi_{0,2} y^2)}{x^2}
\frac{\mathrm{d}\left(\frac{y}{x}\right)}{\frac{y}{x}}.$$
La transformation lin\'eaire diagonale $\left(\phi_{2,0}^{-1}x,(\phi_{2,0}\phi_{0,2})^{-1/2}y\right)$ 
nous permet de supposer que $\phi_{2,0}=\phi_{0,2}=1;$ par cons\'equent
$$\omega=xy\mathrm{d}x+(x^2+\phi_{1,1} xy+y^2)(x\mathrm{d}y-y\mathrm{d}x).$$
Quitte \`a conjuguer $\omega$ par le diff\'eomorphisme $\left(\frac{x}{1
+\phi_{1,1} y}, \frac{y}{1+\phi_{1,1} y}\right)$ le coefficient $\phi_{1,1}$ vaut $0,$ d'o\`u l'\'enonc\'e.
\end{proof}

\medskip
\noindent Remarquons que $\frac{\omega_3}{x^3y}$ est ferm\'ee;
nous en d\'eduisons que $\omega_3$ admet pour int\'egrale premi\`ere
$$\frac{y}{x}\exp\left(\frac{y^2}{2x^2}-\frac{1}{x}\right).$$

\begin{lem}\label{1jetnul3}
{\sl Si le c\^one tangent de $\omega$ en $0$
est r\'eduit \`a une seule droite, alors \`a
conjugaison pr\`es $\mathcal{F}$ est d\'efini par l'une des deux $1$-formes
suivantes
\smallskip
\begin{itemize}
\item[\texttt{1. }] $\omega_1=x^2\mathrm{d}x+y^2(x\mathrm{d}y-y\mathrm{d}x);$

\item[\texttt{2. }] $\omega_2=x^2\mathrm{d}x+(x+y^2)(x\mathrm{d}y-y\mathrm{d}x).$
\end{itemize}}
\end{lem}

\begin{proof}[\sl D\'emonstration]
\`A conjugaison pr\`es, $\mathcal{F}$ est donn\'e par
\begin{align*}
&\omega=x^3\left(\lambda\frac{\mathrm{d}x}{x}+\mathrm{d}\left(\frac{\alpha xy+\beta y^2}{x^2}\right)
\right)+(\phi_{2,0} x^2+\phi_{1,1} xy+\phi_{0,2} y^2)(x\mathrm{d}y-y\mathrm{d}x), && \alpha,\,
\beta,\,\lambda,\,\phi_{i,j}\in\mathbb{C}.
\end{align*}
Le feuilletage $\mathcal{F}$ \'etant quadratique le coefficient $\lambda$ est non nul et nous 
pouvons le supposer \'egal \`a~$1.$ En outre, vu que $0$ est l'unique singularit\'e, nous 
obtenons, en restreignant $\omega$ \`a la droite~$x=0,$ que $\phi_{0,2}\beta=~0.$ 
Puisque $(0:1:0)$ n'est pas singulier $\phi_{0,2}\not=0$
et $\beta=0.$ Quitte \`a conjuguer~$\omega$ par
la transformation lin\'eaire diagonale $(\phi_{0,2}^{1/2}x,y),$ nous nous ramenons \`a $\phi_{0,2}=1.$ Ainsi
$\mathcal{F}$ est d\'ecrit par
$$\omega=x^2\mathrm{d}x+(\alpha x+ \phi_{2,0} x^2+\phi_{1,1} xy+y^2)(x\mathrm{d}y-y\mathrm{d}x).$$
La conjugaison par le diff\'eomorphisme 
$\left(x,y-\frac{\phi_{1,1}}{2}x\right)$ permet d'annuler $\phi_{1,1}.$ Puis
en faisant agir
\begin{align*}
& \left(\frac{x}{1+\phi_{2,0} y},\frac{y}{1+\phi_{2,0} y}\right)  \text{ si }
\alpha=0, && \text{et}&& 
\left(\frac{x}{1-\left(\frac{\phi_{2,0}}{\alpha}\right)x},
\frac{y}{1-\left(\frac{\phi_{2,0}}{\alpha}\right)x}\right) \text{ si } \alpha
\neq 0,
\end{align*}
nous pouvons supposer que $\phi_{2,0}=0.$ Si $\alpha=0$ nous obtenons le premier
mod\`ele $\omega_1.$ Si $\alpha$ est
non nul, alors, en conjugant par la transformation lin\'eaire diagonale 
$(\alpha^3x,\alpha^2y),$ on se ram\`ene \`a $\alpha=1,$ c'est-\`a-dire au
mod\`ele~$\omega_2.$
\end{proof}

\bigskip
\noindent Remarquons que $\frac{\omega_1}{x^4}$ est ferm\'ee et 
a pour int\'egrale premi\`ere 
$\frac{1}{3}\left(\frac{y}{x}\right)^3-\frac{1}{x};$
quant \`a $\omega_2$ elle admet pour int\'egrale premi\`ere
\begin{align*}
\left(2+\frac{1}{x}+2\left(\frac{y}{x}\right)+\left(\frac{y}{x}
\right)^2\right)\exp\left(-\frac{y}{x}\right).
\end{align*}

\noindent Les Lemmes \ref{1jetnul1}, \ref{1jetnul2} et \ref{1jetnul3} impliquent 
l'\'enonc\'e qui suit.

\begin{pro}\label{1jetnul4}
{\sl Soit $\mathcal{F}$ un feuilletage quadratique sur
$\mathbb{P}^2(\mathbb{C})$ dont l'unique singularit\'e est $(0:0:1).$ Si le $1$-jet
de la forme $\omega$ d\'efinissant $\mathcal{F}$ au point singulier est nul alors, 
\`a isomorphisme pr\`es, $\mathcal{F}$
est donn\'e par l'une des trois $1$-formes suivantes
\smallskip
\begin{itemize}
\item[\texttt{1. }] $\omega_1=x^2\mathrm{d}x+y^2(x\mathrm{d}y-y\mathrm{d}x);$

\item[\texttt{2. }] $\omega_2=x^2\mathrm{d}x+(x+y^2)(x\mathrm{d}y-y\mathrm{d}x);$

\item[\texttt{3. }] $\omega_3=xy\mathrm{d}x+(x^2+y^2)(x\mathrm{d}y-y\mathrm{d}x).$
\end{itemize}}
\end{pro}

\medskip

\begin{rem}\label{paslinconj}
Les quatre feuilletages $\mathcal{F}_i$ ne sont pas
lin\'eairement conjugu\'es. En effet 
$\mathcal{F}_4$ est le seul \`a poss\'eder une singularit\'e de type
selle-noeud, il n'est donc pas conjugu\'e aux autres $\mathcal{F}_i.$ La
$1$-forme $\omega_3$ n'est pas lin\'eairement conjugu\'ee
\`a $\omega_1$ et $\omega_2$ car le c\^one tangent en $0$ 
de $\omega_3$ est form\'e de deux droites distinctes, alors que le 
c\^one tangent en $0$ de $\omega_1$ (resp. $\omega_2$) est 
constitu\'e d'une seule droite. Enfin $\mathcal{F}_1$ et 
$\mathcal{F}_2$ ne sont pas lin\'eairement conjugu\'es: le 
premier poss\`ede une int\'egrale premi\`ere rationnelle mais 
pas le second.
\end{rem}

\noindent Le Th\'eor\`eme \ref{feuilquad} r\'esulte des Propositions 
\ref{casnil}, \ref{casselnoeud}, \ref{1jetnul4} et de la Remarque \ref{paslinconj}.\bigskip

\noindent Par contre les feuilletages $\mathcal{F}_2$ et $\mathcal{F}_3$ sont
birationnellement conjugu\'es car ils le sont au feuilletage lin\'eaire d'int\'egrale
premi\`ere $x\mathrm{e}^y;$ leurs feuilles g\'en\'eriques sont 
transcendantes. Dans la carte $X=1$ les feuilletages $\mathcal{F}_1,$ 
$\mathcal{F}_2$ et~$\mathcal{F}_3$ poss\`edent respectivement 
les int\'egrales premi\`eres suivantes
\begin{align*}
& z-\frac{y^3}{3}, &&(2+z+2y+y^2)\mathrm{e}^{-y}, && y\exp\left(\frac{y^2}{2}
-z\right).
\end{align*}
Quitte \`a faire agir les automorphismes polynomiaux de $\mathbb{C}^2$
\begin{align*}
& \left(y,z-\frac{y^3}{3}\right), && \text{resp. }(y,2+z+2y+y^2), && \text{resp. }
\left(y,z-\frac{y^2}{2}\right)
\end{align*}
sur $\mathcal{F}_1,$ resp. $\mathcal{F}_2,$ resp. $\mathcal{F}_3$ nous 
constatons que
\smallskip
\begin{itemize}
\item[\texttt{1. }] $\mathcal{F}_1$ est polynomialement conjugu\'e au feuilletage
d\'ecrit par $z=$ cte;

\item[\texttt{2. }] $\mathcal{F}_2$ et $\mathcal{F}_3$ sont polynomialement 
conjugu\'es au feuilletage donn\'e par les niveaux de $y\mathrm{e}^{-z}
=$~cte.
\end{itemize}
\smallskip

\noindent Nous en d\'eduisons la:
\begin{pro}
{\sl Les feuilletages $\mathcal{F}_1,$ $\mathcal{F}_2$ et $\mathcal{F}_3$
poss\`edent la propri\'et\'e suivante: pour tout point $m$ r\'egulier
la feuille de $\mathcal{F}_i$ passant par $m$ est une courbe 
enti\`ere, {\it i.e.} une courbe param\'etr\'ee par une application holomorphe de 
$\mathbb{C}$ dans $\mathbb{C}^2.$

\noindent En fait $\mathcal{F}_1,$ $\mathcal{F}_2$ et $\mathcal{F}_3$
sont d\'efinis, dans la carte $X=1,$ par des champs de vecteurs polynomiaux
complets.}
\end{pro}

\section{Orbites sous l'action de $\mathrm{PGL}_3(\mathbb{C})$}

\noindent L'ensemble~$\mathscr{F}(2;N)$ est un ouvert de
\textsc{Zariski} dans le projectivis\'e des $1$-formes en les variables $X,$~$Y$ 
et~$Z$ de degr\'e $N+1$ satisfaisant
l'identit\'e d'\textsc{Euler} (il faut en effet que la condition \og sans 
composante commune\fg\hspace{1mm}\'evoqu\'ee dans l'introduction soit satisfaite). En particulier $\mathscr{F}(2;2)$ est un ouvert de 
\textsc{Zariski} dans~$\mathbb{P}^{14}(\mathbb{C}).$ Le groupe des automorphismes
de $\mathbb{P}^2(\mathbb{C})$ agit sur $\mathscr{F}
(2;2);$ l'orbite 
d'un \'el\'ement $\mathcal{F}$ de $\mathscr{F}(2;2)$ sous l'action 
de~$\mathrm{Aut}(\mathbb{P}^2(\mathbb{C}))=\mathrm{PGL}_3(\mathbb{C})$
est not\'ee $\mathcal{O}(\mathcal{F}).$\label{not8} 
Comme cette action est alg\'ebrique
les orbites sont d'adh\'erence (ordinaire) alg\'ebrique dans $\mathbb{P}^{14}
(\mathbb{C}).$ Nous nous int\'eressons dans ce qui suit aux adh\'erences~$\overline{\mathcal{O}(\mathcal{F})}$ des orbites~$\mathcal{O}(\mathcal{F})$ 
dans $\mathscr{F}(2;2)$ sous l'action du 
groupe $\mathrm{Aut}(\mathbb{P}^2(\mathbb{C})).$

\subsection{Isotropies et dimensions des $\mathcal{O}(\mathcal{F}_i)$}\hspace{1mm}

\begin{defi}
Soit $\mathcal{F}$ un feuilletage sur $\mathbb{P}^2(\mathbb{C}).$
Le sous-groupe de $\mathrm{Aut}(\mathbb{P}^2(\mathbb{C}))$
qui pr\'eserve $\mathcal{F}$ s'appelle le {\sl groupe d'isotropie} 
de $\mathcal{F}$ et est not\'e $\mathrm{Iso}(\mathcal{F})$\label{not8d}; 
c'est un groupe alg\'ebrique.
\end{defi}

\noindent La Proposition suivante est de nature \'el\'ementaire, sa 
d\'emonstration est laiss\'ee au lecteur. 

\begin{pro}\label{difdim}
{\sl Les dimensions des $\mathcal{O}(\mathcal{F}_i)$ sont les suivantes
\begin{align*}
& \dim\mathcal{O}(\mathcal{F}_1)=6, && \dim\mathcal{O}(\mathcal{F}_2)=7, &&
\dim\mathcal{O}(\mathcal{F}_3)=7&& \text{et}&&\dim\mathcal{O}(\mathcal{F}_4)=8.
\end{align*}

\noindent Plus pr\'ecis\'ement les groupes $\mathrm{Iso}(\mathcal{F}_i)$
sont donn\'es par
\smallskip
\begin{itemize}
\item[\texttt{1. }]
$\mathrm{Iso}(\mathcal{F}_1)=\{(\beta^3 x:\beta^2y:z+\alpha x)\hspace{1mm}
\vert\hspace{1mm} \alpha\in\mathbb{C},\hspace{1mm}\beta\in\mathbb{C}^*\}$
qui est isomorphe au groupe des transformations
affines de la droite; 

\medskip

\item[\texttt{2. }] $\mathrm{Iso}(\mathcal{F}_2)=\{(x:\alpha x+y:-\alpha(\alpha+2)x-2\alpha y+z)\,\vert\, \alpha\in\mathbb{C}\};$

\medskip

\item[\texttt{3. }] $\mathrm{Iso}(\mathcal{F}_3)=\left\{(x:\pm y:z+\alpha x)
\hspace{1mm}\Big\vert\hspace{1mm}\alpha\in\mathbb{C}\right\};$

\medskip

\item[\texttt{4. }] $\mathrm{Iso}(\mathcal{F}_4)=\{\mathrm{id},\hspace{1mm}(\mathrm{j}x:
\mathrm{j}^2y:z),\hspace{1mm}(\mathrm{j}^2x:\mathrm{j}y:z)\}$
o\`u $\mathrm{j}=\mathrm{e}^{2\mathrm{i}\pi/3}.$
\end{itemize}}
\end{pro}

\subsection{Minoration de la dimension de l'orbite d'un feuilletage quadratique}\hspace{1mm}

\noindent Nous avons vu que la dimension de $\mathcal{O}(\mathcal{F}_1)$
est $6;$ nous montrons dans cette section que c'est la dimension minimale
possible.

\noindent Notons $\chi(\mathbb{P}^2(\mathbb{C}))$\label{not8e} l'alg\`ebre de \textsc{Lie}
des champs de vecteurs holomorphes globaux: $\chi(\mathbb{P}^2(
\mathbb{C}))$ est bien s\^ur l'alg\`ebre de \textsc{Lie} du groupe 
d'automorphismes de $\mathbb{P}^2(\mathbb{C}).$ Soient $\mathcal{F}$ un feuilletage sur
$\mathbb{P}^2(\mathbb{C})$ et $\mathcal{X}$ un 
\'el\'ement de $\chi(\mathbb{P}^2(\mathbb{C}));$ nous dirons que 
$\mathcal{X}$ est une sym\'etrie du feuilletage $\mathcal{F}$ si le
flot $\exp(t\mathcal{X})$ est, pour chaque $t,$ dans le groupe 
d'isotropie $\mathrm{Iso}(\mathcal{F})$ de~$\mathcal{F}.$ Si $\omega$
d\'efinit $\mathcal{F}$ dans une carte affine, $\mathcal{X}$ est une 
sym\'etrie de~$\mathcal{F}$ si et seulement si $L_\mathcal{X}
\omega\wedge\omega=0$ o\`u $L_\mathcal{X}$ d\'esigne la d\'eriv\'ee de \textsc{Lie}
de $\omega$ suivant $\mathcal{X}.$ Remarquons que si $\mathcal{F}$ est de degr\'e
sup\'erieur \`a $2$ et $\mathcal{X}$ une sym\'etrie non triviale de 
$\mathcal{F}$ alors $\mathcal{X}$ n'est pas tangent au feuilletage~$\mathcal{F}:$ si c'\'etait le cas $\mathcal{F}$ serait en effet de 
degr\'e $0$ ou $1.$

\begin{lem}\label{IP}
{\sl Soit $\mathcal{F}$ un \'el\'ement de $\mathscr{F}(2;N)$ d\'efini 
par une $1$-forme $\omega.$
Si $\mathcal{X}$ et
$\mathcal{Y}$ sont deux sym\'etries de~$\mathcal{F}$ ind\'ependantes
sur $\mathbb{C}$ alors $\frac{\omega(\mathcal{X})}{\omega(\mathcal{Y})}$ est 
une int\'egrale premi\`ere rationnelle non constante de 
$\mathcal{F}.$}
\end{lem}

\begin{proof}[\sl D\'emonstration]
D'apr\`es
ce qui pr\'ec\`ede les polyn\^omes~$\omega(\mathcal{X})$ et $\omega(\mathcal{Y})$ sont non identiquement nuls. Les champs 
$\mathcal{X}$ et $\mathcal{Y}$ \'etant ind\'ependants sur $\mathbb{C},$ la fonction rationnelle 
$\frac{\omega(\mathcal{X})}{\omega(\mathcal{Y})}$ est non constante. Rappelons 
que si $\mathcal{X}$ est une sym\'etrie de $\mathcal{F}$ et 
si $\omega(\mathcal{X})\not\equiv 0,$ alors~$\omega(\mathcal{X})$ est un facteur int\'egrant de 
$\omega,$ autrement dit $\frac{\omega}{\omega(\mathcal{X})}$ est ferm\'ee. Puisque le quotient de 
deux facteurs int\'egrants de $\omega$ est une int\'egrale premi\`ere rationnelle de 
$\mathcal{F},$ la fonction $\frac{\omega(\mathcal{X})}{\omega(\mathcal{Y})}$ est une int\'egrale premi\`ere
rationnelle non constante de $\mathcal{F}.$
\end{proof}

\noindent Si $\mathrm{G}$ d\'esigne un groupe alg\'ebrique, nous notons 
$\mathfrak{g}$\label{not9} son alg\`ebre de \textsc{Lie}.

\begin{pro}\label{dimmin}
{\sl Soient $N$ un entier sup\'erieur ou \'egal \`a $2$ et $\mathcal{F}$ un \'el\'ement 
de $\mathscr{F}(2;N).$ La dimension de~$\mathcal{O}(\mathcal{F})$ 
est minor\'ee par $6.$ }
\end{pro}

\begin{proof}[\sl D\'emonstration] 
Comme $\dim\chi(\mathbb{P}^2(\mathbb{C}))=8$ le Lemme \ref{IP} assure que si $\dim\mathcal{O}(\mathcal{F})\leq 6,$ le feuilletage~$\mathcal{F}$ 
admet une int\'egrale premi\`ere
rationnelle. Dans \cite{CM} \textsc{Cerveau} et \textsc{Mattei} montrent, en utilisant 
le th\'eor\`eme de \textsc{L{\"u}roth}, qu'il existe
une fonction rationnelle non constante $\mathrm{f},$ d\'efinie \`a~$\mathrm{PGL}_2(\mathbb{C})$
pr\`es, telle que l'ensemble des int\'egrales premi\`eres rationnelles de 
$\mathcal{F}$ soit isomorphe \`a $\mathbb{C}(\mathrm{f});$
on dit alors que $\mathrm{f}$ est minimale. Choisissons~$\mathrm{f}$ minimale pour $\mathcal{F}.$ 
Remarquons que si $\varphi$ appartient \`a $\mathrm{Iso}(\mathcal{F})$ alors~$\mathrm{f}\circ\varphi$ est encore
une int\'egrale premi\`ere de~$\mathcal{F}$ et que, de plus, elle est minimale;~$\mathrm{f}\circ
\varphi$ s'\'ecrit donc~$\tau_\varphi(\mathrm{f})$ avec $\tau_\varphi$ dans $\mathrm{PGL}_2(\mathbb{C}).$
D\'esignons par $\tau$ le morphisme d\'efini par
\begin{align*}
&\tau\colon\mathrm{Iso}(\mathcal{F})\to\mathrm{PGL}_2(\mathbb{C}), && \varphi\mapsto \tau_\varphi.
\end{align*} 
Le noyau de $\tau$ est un sous-groupe de $\mathrm{Iso}(\mathcal{F})$ n\'ecessairement discret.
En effet si 
$\dim\ker \tau\geq 1,$ alors~$\ker \tau$ contient un flot dont le g\'en\'erateur infinit\'esimal est 
tangent \`a $\mathcal{F}$ ce qui est impossible. Par suite 
$$\mathrm{D}_{id}\tau\colon\mathfrak{iso}(\mathcal{F})\to\mathfrak{sl}_2(\mathbb{C})$$
est injective; ceci implique l'in\'egalit\'e 
$$\dim\mathrm{Iso}(\mathcal{F})=\dim\mathfrak{iso}(\mathcal{F})\leq\dim\mathfrak{sl}_2(\mathbb{C})=3.$$
Supposons que $\dim\mathrm{Iso}(\mathcal{F})=3$ alors $\mathfrak{iso}
(\mathcal{F})\simeq\mathfrak{sl}_2
(\mathbb{C}).$ Soit $\mathrm{G}$ la composante neutre de 
$\mathrm{Iso}(\mathcal{F}).$ \`A partir de~$\mathfrak{iso}(\mathcal{F})
=\mathfrak{sl}_2(\mathbb{C})$
nous obtenons \`a isomorphisme pr\`es 
\begin{align*}
& \mathrm{G}=\mathrm{SL}_2(\mathbb{C})&&\text{ou} &&\mathrm{G}=
\mathrm{PGL}_2(\mathbb{C}).
\end{align*}
Dans chacune de ces \'eventualit\'es nous h\'eritons d'une action
$\rho$ de $\mathrm{SL}_2(\mathbb{C})$ sur $\mathbb{P}^2(\mathbb{C}).$ 

\noindent Si $\rho$ est irr\'eductible $\rho$ s'identifie au projectivis\'e de l'action
naturelle $\widetilde{\rho}$ de $\mathrm{SL}_2(\mathbb{C})$ sur l'espace des 
formes quadratiques en deux variables. Cette action n'a pas de feuilletage invariant 
de degr\'e sup\'erieur ou \'egal \`a $2;$ en effet si $\mathcal{F}$ \'etait invariant, son lieu singulier,
qui n'est pas vide, le serait aussi: contradiction avec le fait que $\widetilde{\rho}$ 
ne laisse invariant aucun ensemble fini. 

\noindent Si $\rho$ est r\'eductible, elle se d\'ecompose, au niveau de $\mathbb{C}^3,$
en une action sur $\mathbb{C}^2$ et une sur $\mathbb{C}$ qui est n\'ecesseraiment triviale.
Par suite dans une carte affine ad-hoc $\rho$ co\"incide avec l'action lin\'eaire standard
$$
\left(\left(
\begin{array}{cc}
\alpha & \beta\\
\gamma & \delta
\end{array}
\right),(x,y)\right)\mapsto(\alpha x+\beta y,\gamma x+\delta y).$$
Cette derni\`ere action pr\'eserve un feuilletage: celui induit par le champ radial 
$\mathcal{R}=x\frac{\partial}{\partial x}+~y\frac{\partial}{\partial y}$\label{not10}, qui est 
de degr\'e~$0,$ et c'est le seul. En effet soit $\widetilde{\mathcal{F}}$ un feuilletage non radial
invariant par 
$$(x,y)\mapsto(\alpha x+\beta y,\gamma x+\delta y);$$
il existe alors un point g\'en\'erique $m$ tel que la tangente \`a la feuille
$\widetilde{\mathcal{L}}_m$ passant par $m$ ne soit pas la droite passant par 
l'origine et $m.$ Nous pouvons supposer que $m=(1,0);$ le groupe 
$$\{(x+\beta y,\delta y)\hspace{1mm}\vert\hspace{1mm}\beta\in\mathbb{C},\hspace{1mm}
\delta\in\mathbb{C}^*\}$$
fixe $m,$ doit fixer $\widetilde{\mathcal{F}}$ donc doit fixer la tangente \`a $\widetilde{\mathcal{L}}_m$
en $m$ ce qui est impossible. 
Nous en d\'eduisons que~$\dim\mathrm{Iso}(\mathcal{F})\not=3.$ Il s'en suit que~$\dim\mathrm{Iso}(\mathcal{F})\leq 2$ et $\dim\mathcal{O}(\mathcal{F})\geq~6.$
\end{proof} 

\begin{cor}
{\sl Le feuilletage $\mathcal{F}_1$ r\'ealise la dimension minimale des 
orbites en degr\'e $2.$}
\end{cor}

\noindent Nous pouvons nous demander si $\mathcal{F}_1$ est l'unique
\'el\'ement de $\mathscr{F}(2;2)$ \`a satisfaire cette propri\'et\'e; nous allons
voir dans le paragraphe suivant que ce n'est pas le cas.

\subsection{Description des feuilletages quadratiques tels que
$\dim\mathcal{O}(\mathcal{F})=6$}\hspace{1mm}

\begin{pro}\label{nonab}
{\sl Soit $\mathcal{F}$ un feuilletage de degr\'e $N$ sur $\mathbb{P}^2
(\mathbb{C}).$ Si $\dim\mathfrak{iso}(\mathcal{F})=2,$
alors $\mathfrak{iso}(\mathcal{F})$ est isomorphe \`a l'alg\`ebre
du groupe des transformations affines de la droite.}
\end{pro}

\noindent 

\begin{proof}[\sl D\'emonstration]
D'apr\`es la classification des alg\`ebres de \textsc{Lie} de dimension $2$
(\emph{voir} \cite{FH}) il suffit de montrer que $\mathfrak{iso}(\mathcal{F})$ 
n'est pas ab\'elienne. Pour ce 
faire raisonnons par l'absurde, {\it i.e.} supposons que~$\mathfrak{iso}(\mathcal{F})$ soit ab\'elienne. Soient $\mathcal{X}$ et $\mathcal{Y}$ deux g\'en\'erateurs de $\mathfrak{iso}(\mathcal{F}).$ La triangulation des 
alg\`ebres r\'esolubles de matrices assure l'existence d'une droite $\mathcal{D}$ invariante
par $\mathcal{X}$ et~$\mathcal{Y}.$ Nous pouvons alors nous ramener \`a la 
situation suivante: $\mathcal{D}$
est la droite \`a l'infini et dans la carte~$\mathbb{C}^2$ les 
champs $\mathcal{X}$ et~$\mathcal{Y}$ sont affines.
Dans cette m\^eme carte $\mathcal{F}$ est d\'efini par un champ polynomial
\`a singularit\'e isol\'ee que nous noterons~$\mathcal{Z}.$ Soit~$\mathrm{V}$ un \'el\'ement de~$\mathfrak{iso}
(\mathcal{F}).$ Puisque
le flot de~$\mathrm{V}$ laisse $\mathcal{F}$ invariant, en
tout point non singulier~$[\mathrm{V},\mathcal{Z}]$ est un multiple de $\mathcal{Z};$ d'apr\`es \cite{S} il existe une 
fonction holomorphe $g$ sur~$\mathbb{C}^2$ telle que~$[\mathrm{V},\mathcal{Z}]=~g \mathcal{Z}.$
Les champs~$\mathrm{V}$ et $\mathcal{Z}$ \'etant polynomiaux, $[\mathrm{V},\mathcal{Z}]$
l'est aussi; par suite $g$ est rationnelle et holomorphe sur~$\mathbb{C}^2$ 
donc polynomiale. Le champ $\mathrm{V}$ \'etant affine nous 
avons $\deg [\mathrm{V},\mathcal{Z}]\leq \deg \mathcal{Z};$ la relation~$[\mathrm{V},
\mathcal{Z}]=~g \mathcal{Z}$
entra\^ine que $g$ est constante. Finalement nous obtenons 
\begin{align*}
&[\mathcal{X},\mathcal{Z}]=\lambda \mathcal{Z}, &&[\mathcal{Y},\mathcal{Z}]=\eta \mathcal{Z}, &&\lambda, \hspace{1mm}\eta\in 
\mathbb{C}.
\end{align*}
\`A combinaison lin\'eaire et permutation pr\`es de $\mathcal{X}$ et $\mathcal{Y}$
nous pouvons nous ramener \`a  
\begin{align*}
&[\mathcal{X},\mathcal{Z}]=0, &&[\mathcal{Y},\mathcal{Z}]=\eta \mathcal{Z}, &&\eta\in 
\mathbb{C}.
\end{align*} 
Commen\c{c}ons par supposer que $\mathcal{X}$ et $\mathcal{Y}$ sont g\'en\'eriquement 
transverses. En un point g\'en\'erique il existe des coordonn\'ees locales 
telles que 
\begin{align*}
& \mathcal{X}=\frac{\partial}{\partial x}, &&  \mathcal{Y}=\frac{\partial}{\partial y}. 
\end{align*}
Dans ce syst\`eme de coordonn\'ees $\mathcal{Z}$ est de la forme
\begin{align*}
&\alpha(x,y) \frac{\partial}{\partial x}+\beta(x,y)\frac{\partial}{\partial y}, && 
\alpha,\hspace{1mm}\beta\in\mathbb{C}\{x,y\}.
\end{align*}
L'\'egalit\'e $[\mathcal{X},\mathcal{Z}]=0$ implique que $\alpha$ et $\beta$ ne d\'ependent 
pas de $x.$ Alors $[\mathcal{Y},\mathcal{Z}]=\eta \mathcal{Y}$ entra\^ine que 
$\alpha'(y)=\eta\alpha(y)$ et $\beta'(y)=\eta\beta(y).$
On obtient ainsi 
\begin{align*}
&\mathcal{Z}=e^{\eta y}\left(\gamma\frac{\partial}{\partial x}+\delta\frac{
\partial}{\partial y}\right), && \gamma,\hspace{1mm}\delta\in\mathbb{C}, 
\end{align*}
qui se r\'e\'ecrit $\mathcal{Z}=e^{\eta y}(\gamma\mathcal{X}+\delta\mathcal{Y}).$
Sur un ouvert $\mathcal{Z}$ et $\gamma\mathcal{X}+\delta\mathcal{Y}$ sont donc parall\`eles; nous
en d\'eduisons que c'est le cas partout. Ainsi le feuilletage d\'efini par $\mathcal{Z}$ est
lin\'eaire: contradiction.

\noindent Supposons que $\mathcal{X}$ et $\mathcal{Y}$ soient partout colin\'eaires; par cons\'equent en un 
point g\'en\'erique nous pouvons \'ecrire 
\begin{align*}
& \mathcal{X}=\frac{\partial}{\partial x}, && \mathcal{Y}=\delta(y)\frac{\partial}{\partial x}, &&
\delta\in\mathbb{C}\{y\}.
\end{align*}
Notons que $\mathcal{X}$ et $\mathcal{Y}$ n'\'etant pas $\mathbb{C}$-colin\'eaires, 
$\delta$ est non constante.
\`A partir de $[\mathcal{X},\mathcal{Z}]=0,$ nous obtenons 
\begin{align*}
& \mathcal{Z}=\alpha(y)\frac{\partial}{\partial x}+\beta(y)\frac{\partial}{\partial y}, && 
\alpha,\hspace{1mm}\beta\in\mathbb{C}\{y\};
\end{align*}
remarquons que $\beta\not\equiv 0$ sinon $\mathcal{Z}$ serait colin\'eaire \`a $\mathcal{X}$
et d\'efinirait un feuilletage lin\'eaire. Par ailleurs~$[\mathcal{Y},\mathcal{Z}]=~\eta \mathcal{Z}$ conduit \`a $\eta\alpha(y)=0$ et $\beta(y)(\eta+\delta'(y))=0,$
soit \`a $\alpha\equiv 0$ et $\delta'(y)=-\eta\not=0.$
Autrement dit 
\begin{align*}
& \mathcal{X}=\frac{\partial}{\partial x}, && \mathcal{Y}=(-\eta y+\varepsilon)\frac{\partial}{\partial x},
&& \mathcal{Z}=\beta(y)\frac{\partial}{\partial y}, &&\eta\in\mathbb{C}^*,\hspace{1mm}
\varepsilon\in\mathbb{C},\hspace{1mm}\beta\in\mathbb{C}\{y\}.
\end{align*}
Par suite $\mathcal{F}$ a une int\'egrale premi\`ere rationnelle $\mathrm{f}$ (Lemme
\ref{IP}) qui, dans les coor\-donn\'ees locales~$(x,y)$ ne d\'epend pas de $y$ et s'\'ecrit 
donc~$\mathrm{f}(x).$ Le flot de $\mathcal{Y}$ est donn\'e par $(x+t\delta(y),y);$
puisque $\mathcal{Y}$ est une sym\'etrie de~$\mathcal{F},$ la compos\'ee~$\mathrm{f}(x+
t\delta(y))$ ne d\'epend pas de $y,$ {\it i.e.} $\delta$ est une constante
ce qui, comme nous l'avons d\'ej\`a vu, est impossible.
\end{proof}

\begin{lem}\label{lieusing}
{\sl Soit $\mathcal{F}$ un feuilletage quadratique sur le plan projectif
complexe. 
Supposons que l'alg\`ebre de \textsc{Lie} $\mathfrak{iso}(\mathcal{F})$ 
soit de type affine, {\it i.e.} engendr\'ee par deux champs de
vecteurs $\mathcal{X}$ et $\mathcal{Y}$ tels que $[\mathcal{X},\mathcal{Y}]=
\mathcal{Y}.$ Il existe une droite dans $\mathbb{P}^2(\mathbb{C})$
invariante par les champs $\mathcal{X}$ et $\mathcal{Y}$ et telle que tout point
singulier de $\mathcal{F}$ appartienne \`a cette droite.}
\end{lem}

\begin{proof}[\sl D\'emonstration]
Puisque $[\mathcal{X},\mathcal{Y}]=\mathcal{Y},$  le th\'eor\`eme d'\textsc{Engel}
assure l'existence d'une droite~$\mathcal{D}$ de $\mathbb{P}^2(\mathbb{C}),$ 
que nous pouvons supposer \^etre la droite $Z=0,$ invariante par $\mathcal{X}$ et~$\mathcal{Y}.$
Pla\c{c}ons-nous dans la car\-te~$Z=1;$ en reprenant un argument de la 
d\'emonstration de la Proposition \ref{nonab} nous avons
\begin{align*}
&[\mathcal{X},\mathcal{Z}]=\lambda \mathcal{Z}, &&[\mathcal{Y},\mathcal{Z}]=\eta \mathcal{Z}, &&\lambda,\hspace{1mm}\eta\in \mathbb{C}.
\end{align*} 
L'identit\'e de \textsc{Jacobi} 
$$[\mathcal{X},[\mathcal{Y},\mathcal{Z}]]+[\mathcal{Y},[\mathcal{Z},\mathcal{X}]]+[\mathcal{Z},[\mathcal{X},\mathcal{Y}]]=0$$
entra\^ine $\eta=0,$ c'est-\`a-dire $[\mathcal{Y},\mathcal{Z}]=0.$

\noindent Si toutes les singularit\'es de $\mathcal{F}$ sont sur la droite \`a l'infini $\mathcal{D}$ le lemme est d\'emontr\'e.

\noindent Supposons qu'il existe un point singulier $m$ \`a distance finie pour le feuilletage $\mathcal{F};$ 
\'evidemment ce point est singulier pour $\mathcal{X}$ et $\mathcal{Y}.$ Nous pouvons nous ramener \`a~$m=(0,0)$ et 
\begin{align*}
& \mathcal{X}=-(1+\alpha)x\frac{\partial}{\partial x}-\alpha y\frac{\partial}{\partial y}, &&
\mathcal{Y}=y\frac{\partial}{\partial x}.
\end{align*}

\noindent Quant \`a $\mathcal{Z}$ nous l'\'ecrivons sous la forme $\mathcal{Z}_1+\mathcal{Z}_2+\phi \mathcal{R}$ o\`u 
$\mathcal{Z}_1$ (resp. $\mathcal{Z}_2$) 
d\'esigne un champ lin\'eaire (resp. quadratique) et $\mathcal{R}$ le champ radial. L'\'egalit\'e
$[\mathcal{Y},\mathcal{Z}]=0$ conduit \`a: $[\mathcal{Y},\mathcal{Z}_1]=[\mathcal{Y},\mathcal{Z}_2]=0$ et~$\phi$ est du type $\varepsilon y^2.$
\`A partir de $[\mathcal{Y},\mathcal{Z}_1]=0$ et $[\mathcal{Y},\mathcal{Z}_2]=0$ nous obtenons
\begin{align*}
&\mathcal{Z}_1=(\beta x+\gamma y)\frac{\partial}{\partial x}+\beta y\frac{\partial}{\partial y}, &&
\mathcal{Z}_2=(\delta y^2+\kappa xy)\frac{\partial}{\partial x}+\kappa y^2\frac{\partial}{\partial y}, 
&& \beta,\hspace{1mm} \gamma, \hspace{1mm} \kappa, \hspace{1mm}
\delta\in\mathbb{C}.
\end{align*}
Nous constatons que si $\beta$ est nul, $\mathcal{Z}$ est divisible par 
$y,$ donc n'est pas associ\'e \`a un feuilletage quadratique. Il 
s'en suit que $\beta$ est non nul ce qui implique que $\lambda=0.$
En \'ecrivant explicitement~$[\mathcal{X},\mathcal{Z}]=0$ nous
avons
$$\gamma=\alpha\varepsilon=\alpha\kappa=\delta(\alpha-1)=0.$$
Si $\delta$ est nul, $\mathcal{Z}$ est colin\'eaire au champ radial 
ce qui ne 
peut arriver. Ainsi $\alpha=1$ et~$\varepsilon=\kappa=0;$ d'o\`u
$$\mathcal{Z}=\beta\mathcal{R}+\delta y^2\frac{\partial}{\partial x}.$$
Nous pouvons \'evidemment normaliser les coefficients $\beta$ et $\delta$
\`a $1.$ Le feuilletage $\mathcal{F}$ est d\'efini par la $1$-forme
$y\mathrm{d}x-x\mathrm{d}y-y^2\mathrm{d}y$ qui poss\`ede l'int\'egrale premi\`ere $\frac{x}{y}-y;$
ses points singuliers sont~$(0:0:1)$ et $(1:0:0)$ et la 
droite $y=0$ satisfait l'\'enonc\'e.
\end{proof}

\noindent Les feuilletages quadratiques dont l'orbite est de dimension maximale $6$
sont classifi\'es par la:

\begin{pro}\label{dim6}
{\sl Soit $\mathcal{F}$ un \'el\'ement de $\mathscr{F}(2;2).$
Supposons
que l'alg\`ebre de \textsc{Lie} $\mathfrak{iso}(\mathcal{F})$ soit de 
dimension $2;$ alors~$\mathcal{F}$ est, \`a action d'un automorphisme 
de~$\mathbb{P}^2(\mathbb{C})$ pr\`es, d\'efini par l'une des $1$-formes
suivantes
\begin{itemize}
\item[\texttt{1. }] $\omega_1=x^2\mathrm{d}x+y^2(x\mathrm{d}y-y\mathrm{d}x);$

\item[\texttt{2. }] $\omega_5=x^2\mathrm{d}y+y^2(x\mathrm{d}y-y\mathrm{d}x).$
\end{itemize}

\noindent Dit autrement les orbites associ\'ees sont les seules orbites de dimension $6.$
Elles sont ferm\'ees dans $\mathscr{F}(2;2).$

\noindent De plus nous avons
\begin{align*}
&\mathrm{Iso}(\mathcal{F}_1)=\{(\alpha^3 x,\alpha^2y),\hspace{1mm}
\left(\frac{x}{1+\beta x},\frac{y}{1+\beta x}\right)\hspace{1mm}
\Big\vert\hspace{1mm} \alpha\in\mathbb{C}^*,\hspace{1mm}\beta\in\mathbb{C}\},\\
& \mathrm{Iso}(\mathcal{F}_5)=\left\{(\alpha^2 x,\alpha y),\hspace{1mm}
 \left(\frac{x}{1+\beta y},\frac{y}{1+\beta y}\right) \hspace{1mm}
\Big\vert\hspace{1mm} \alpha\in\mathbb{C}^*,
\hspace{1mm}\beta\in\mathbb{C}\right\}; 
\end{align*}
ces deux groupes ne sont pas conjugu\'es.}
\end{pro}

\begin{rem}
Le feuilletage $\mathcal{F}_5$ d\'ecrit par $\omega_5$ intervient 
d\'ej\`a dans la d\'emonstration du Lemme \ref{lieusing}.
\end{rem}

\begin{proof}[\sl D\'emonstration]
D'apr\`es la Proposition \ref{nonab} et le Lemme \ref{lieusing} il existe 
deux sym\'etries~$\mathcal{X}$ et $\mathcal{Y}$ de~$\mathcal{F}$
telles que $[\mathcal{X},\mathcal{Y}]=\mathcal{Y}.$ Ces deux sym\'etries
pr\'eservent une droite $\mathcal{D}$ que nous supposerons \^etre~$y=0.$ 
De plus $\mathrm{Sing}(\mathcal{F})$ est contenu 
dans $\mathcal{D}.$ Si $\mathcal{F}$ a un unique point singulier
nous savons d'apr\`es le Th\'eor\`eme \ref{feuilquad} et la Proposition 
\ref{difdim} que $\mathcal{F}$ est conjugu\'e \`a $\mathcal{F}_1.$
Nous nous ramenons donc au cas o\`u $\mathrm{Sing}(\mathcal{F})$ 
contient $(0:0:1)$ et $(1:0:0).$ Notons que $\mathcal{X}$ et $\mathcal{Y}$ 
sont singuliers en ces deux points; nous en d\'eduisons en utilisant l'identit\'e~$[\mathcal{X},\mathcal{Y}]=\mathcal{Y}$ que 
\begin{align*}
& \mathcal{X}=(\lambda x+\varepsilon y)\frac{\partial}{\partial x}+\eta y\frac{\partial}{\partial 
y}+\alpha y\mathcal{R},
&& \mathcal{Y}=\gamma y\frac{\partial}{\partial x}+\beta y\mathcal{R}
&&\text{avec} &&(\eta-\lambda)\gamma=\gamma && \text{et} && \eta\beta=\beta
\end{align*}
avec $\lambda,$ $\varepsilon,$ $\eta,$ $\alpha,$ $\beta$ et 
$\gamma$ dans $\mathbb{C}.$ Consid\'erons un 
champ de vecteurs polynomial $\mathcal{Z}$ d\'efinissant $\mathcal{F}$ 
dans la carte~$Z=1;$ \'ecrivons-le sous la forme
$$\mathcal{Z}=\mathcal{Z}_1+\mathcal{Z}_2+\phi\mathcal{R}$$
o\`u $\mathcal{Z}_i$ est un champ de vecteurs homog\`ene
de degr\'e $i$ et $\phi$ une forme quadratique. Nous pouvons 
supposer que la droite \`a l'infini n'est pas invariante par $\mathcal{F}$
ce qui s'interpr\`ete comme suit: $\phi$ n'est pas identiquement 
nulle. Puisque $\mathcal{F}$ est singulier en $(1:0:0),$ la forme 
quadratique $\phi$ est du type suivant $$\phi=\phi_{1,1}xy+\phi_{0,2}y^2.$$

\noindent En reprenant un argument \'evoqu\'e pr\'ec\'edemment
le fait que $\mathcal{X}$ et $\mathcal{Y}$ soient des sym\'etries
de $\mathcal{F}$ se traduit par
\begin{equation}\label{sym1}
[\mathcal{X},\mathcal{Z}]=(a_1+b_1x+c_1y)\mathcal{Z}
\end{equation}
\begin{equation}\label{sym2}
[\mathcal{Y},\mathcal{Z}]=(a_2+b_2x+c_2y)\mathcal{Z}
\end{equation}
les $a_i,$ $b_i$ et $c_i$ d\'esignant des complexes.

\noindent De nouveau on utilise l'identit\'e de \textsc{Jacobi} $[\mathcal{X},[\mathcal{Y},
\mathcal{Z}]]+[\mathcal{Y},[\mathcal{Z},\mathcal{X}]]+[\mathcal{Z},[
\mathcal{X},\mathcal{Y}]]=0$ qui implique ici
$$\mathcal{X}(a_2+b_2x+c_2y)-\mathcal{Y}(a_1+b_1x+c_1y)=a_2+b_2x+c_2y.$$

\noindent Les champs $\mathcal{X}$ et $\mathcal{Y}$ s'annulant
en $0,$ nous obtenons en particulier que $a_2$ est nul puis que 
$$\alpha(b_2x+c_2y)=\beta(b_1x+c_1y).$$ L'\'etude des
termes de plus haut degr\'e dans (\ref{sym1}) conduit \`a 
$[\alpha y\mathcal{R},\phi\mathcal{R}]=(b_1x+c_1y)\phi 
\mathcal{R};$ puisque $\phi$ est non nulle nous en d\'eduisons que
$$\alpha y=b_1x+c_1y.$$
De m\^eme en consid\'erant les termes de plus haut degr\'e 
dans (\ref{sym2}) nous obtenons $$\beta y=b_2x+c_2y.$$
Autrement dit 
\begin{align*}
& b_1=b_2=0, && \alpha=c_1 && \text{et} &&\beta=c_2.
\end{align*}

\noindent La preuve se s\'epare en deux suivant la nullit\'e ou non 
de $\gamma.$

\noindent Dans un premier temps supposons que $\gamma$ soit nul;
nous pouvons alors nous ramener \`a $\beta=1$ et, quitte \`a 
soustraire $\alpha\mathcal{Y}$ \`a $\mathcal{X},$ \`a $\alpha=0.$ 
En r\'e\'ecrivant la relation $[\mathcal{X},\mathcal{Y}]=\mathcal{Y}$
nous constatons que~$\eta=1,$~{\it i.e.} 
\begin{align*}
&\mathcal{X}=(\lambda x+\varepsilon y)\frac{\partial}{\partial x}+y
\frac{\partial}{\partial y}, && \mathcal{Y}=y\mathcal{R}.
\end{align*}

\noindent Les \'egalit\'es (\ref{sym1}) et (\ref{sym2}) sont d\'esormais
du type 
\begin{align*}
&[\mathcal{X},\mathcal{Z}]=a_1\mathcal{Z} &&\text{et} && [\mathcal{Y},
\mathcal{Z}]=y\mathcal{Z}.
\end{align*}

\noindent Cette derni\`ere entra\^ine que
$$-\mathcal{Z}_1(y)\mathcal{R}+y\mathcal{Z}_2-\mathcal{Z}_2(y)\mathcal{R}
+y\phi\mathcal{R}=y(\mathcal{Z}_1+\mathcal{Z}_2+\phi\mathcal{R})$$
d'o\`u $-\mathcal{Z}_1(y)\mathcal{R}=y\mathcal{Z}_1$ et $\mathcal{Z}_2
(y)=0.$ Il s'en suit que 
\begin{align*}
&\mathcal{Z}=Q\frac{\partial}{\partial x}+\phi\mathcal{R} && \text{avec} &&Q=q_{2,0}x^2+q_{1,1}xy+
q_{0,2}y^2,
\end{align*}
les formes quadratiques $Q$ et $\phi$ \'etant non nulles et sans facteur commun;
ce qui implique en particulier la non nullit\'e de $q_{2,0}.$

\noindent L'identit\'e $[\mathcal{X},\mathcal{Z}]=a_1\mathcal{Z}$ conduit 
au syst\`eme
\begin{align*}
(\diamond)\left\{\begin{array}{ll}
(\lambda x+\varepsilon y)\frac{\partial Q}{\partial x}+y\frac{\partial Q}{
\partial y}=(a_1+\lambda)Q\\
(\lambda x+\varepsilon y)\frac{\partial\phi}{\partial x}+y\frac{\partial\phi}{
\partial y}=a_1\phi
\end{array}
\right.
\end{align*}

\noindent Montrons que nous pouvons supposer que $\varepsilon=0.$
Remarquons que si $\lambda$ vaut $1,$ la d\'ecomposition de 
\textsc{Jordan} de l'op\'erateur de d\'erivation $\mathcal{X}$ (pour lequel 
$Q$ et $q$ sont vecteurs propres) implique
que 
\begin{align*}
&\varepsilon y\frac{\partial Q}{\partial x}=0 &&\text{et}
&& \varepsilon y\frac{\partial\phi}{\partial x}=0.
\end{align*}

\noindent Puisque $Q$ et $\phi$ sont non nulles et sans facteur commun, 
$\varepsilon$ est n\'ecessairement nul. Enfin si $\lambda$
est diff\'erent de $1,$ nous invoquons la diagonalisation pour nous
ramener \`a $\varepsilon=0;$ notons que cette diagonalisation 
n'alt\`ere pas la forme de $\phi,$ {\it i.e.} $\phi$ est toujours
du type $\phi_{1,1} xy+\phi_{0,2} y^2.$ 

\noindent Le syst\`eme $(\diamond)$ se r\'e\'ecrit alors 
\begin{align*}
(\lambda+1-a_1)\phi_{1,1}=(2-a_1)\phi_{0,2}=(\lambda-a_1)q_{2,0}=(1-a_1)q_{1,1}
=(2-\lambda-a_1)q_{0,2}=0.
\end{align*}

\noindent Comme le coefficient $q_{2,0}$ est non nul on a $\lambda=a_1$ et $\phi_{1,1}=0.$
La forme quadratique $\phi$ \'etant non triviale $\phi_{0,2}$ est non nul et on obtient
$\lambda=a_1=2$ et $q_{1,1}=q_{2,0}=0.$

\noindent Par cons\'equent 
$$\mathcal{Z}=q_{2,0}x^2\frac{\partial}{\partial x}+\phi_{0,2}y^2\mathcal{R}$$
et nous pouvons \'evidemment normaliser les coefficients $q_{2,0}$
et $\phi_{0,2}$ \`a $1.$ Nous obtenons alors le feuilletage $\mathcal{F}_5$ 
d\'efini par $\omega_5.$

\noindent Maintenant nous \'etudions la deuxi\`eme \'eventualit\'e o\`u $\gamma\not=0;$ on peut supposer que $\gamma=1.$ Dans ce cas on a~$\eta=\lambda+1.$ Quitte \`a remplacer $\mathcal{X}$ 
par $\mathcal{X}-\varepsilon\mathcal{Y}$ on se ram\`ene~\`a
\begin{align*}
& \mathcal{X}=\lambda x\frac{\partial}{\partial x}+(\lambda+1)y\frac{
\partial}{\partial y}+\alpha y\mathcal{R}
&& \text{et} && \mathcal{Y}=y\frac{\partial}{\partial x}+\beta y \mathcal{R}.
\end{align*}
Les \'egalit\'es $[\mathcal{X},\mathcal{Z}]=(a_1+\alpha y)\mathcal{Z}$ 
et $[\mathcal{Y},\mathcal{Z}]=\beta y\mathcal{Z}$ ajout\'ees \`a l'identit\'e de \textsc{Jacobi}
nous conduisent \`a~$\lambda\beta=0.$

\noindent Si $\beta=0$ nous remarquons que le champ $\mathcal{Y}=
y\frac{\partial}{\partial x}$ est conjugu\'e via un automorphisme de 
$\mathbb{P}^2(\mathbb{C})$ \`a~$y\mathcal{R}.$ En permutant le 
r\^ole des points singuliers nous nous ramenons au premier cas 
d\'ej\`a trait\'e.

\noindent Lorsque $\lambda=0$ on peut supposer par conjugaison que
$\beta$ vaut $1$ ce qui nous conduit \`a  
\begin{align*}
&\mathcal{X}=-\alpha y\frac{\partial}{\partial x}+y\frac{\partial}{\partial y}
&& \text{et} &&\mathcal{Y}=y\frac{\partial}{\partial x}+y\mathcal{R},
\end{align*}

\noindent quitte \`a soustraire $\alpha\mathcal{Y}$ \`a $\mathcal{X}.$

\noindent Via une transformation lin\'eaire de $\mathbb{C}^2$ qui 
n'affecte ni le champ $\mathcal{Y}$ ni le fait que $\phi$ soit non nul,
nous pouvons diagonaliser $\mathcal{X}$
$$\mathcal{X}=y\frac{\partial}{\partial y};$$
on v\'erifie alors que $[\mathcal{X},\mathcal{Z}]=
a_1\mathcal{Z}.$ Avec les notations habituelles nous avons
$\phi_{1,1}(1-a_1)=\phi_{0,2}(2-a_1)=0$
de sorte que $a_1$ vaut $1$ ou $2.$ 

\noindent Si $a_1=1$ le coefficient $\phi_{0,2}$ est nul et on constate
que $\mathcal{Z}$ est divisible par $y$ et ne d\'efinit
donc pas un feuilletage quadratique.

\noindent Si $a_1=2,$ alors $\phi_{1,1}=0$ et, toujours en invoquant
le fait que $[\mathcal{X},\mathcal{Z}]=a_1\mathcal{Z},$ on v\'erifie 
encore que $\mathcal{Z}$ est divisible par $y.$ 
L'\'eventualit\'e $\gamma\not=0$ n'arrive pas et la 
proposition est d\'emontr\'ee.
\end{proof}

\subsection{Adh\'erences d'orbites}\hspace{1mm}

\noindent Nous commen\c{c}ons l'\'etude des adh\'erences 
d'orbites par celles des $\mathcal{O}(\mathcal{F}_i)$ 
dans $\mathscr{F}(2;2)$ pour $i=1,$ $2,$ $3$ et $5;$ celle
de $\mathcal{O}(\mathcal{F}_4)$  r\'esultera d'un principe
plus g\'en\'eral. 

\noindent La d\'efinition suivante nous sera utile.

\begin{defi}
Soient $\mathcal{F}$ et $\mathcal{F}'$ deux feuilletages. 
On dit que $\mathcal{F}$ {\sl d\'eg\'en\`ere} sur $\mathcal{F}'$ 
si l'adh\'erence $\overline{\mathcal{O}(\mathcal{F})}$ (dans~$\mathscr{F}(2;2)$) de $\mathcal{O}(\mathcal{F})$ contient~$\mathcal{O}(\mathcal{F}')$ et 
$\mathcal{O}(\mathcal{F})\not=\mathcal{O}(\mathcal{F}').$
\end{defi}

\begin{rem}
Si $\mathcal{F}$ d\'eg\'en\`ere sur $\mathcal{F}'$ nous avons 
l'in\'egalit\'e $\dim\mathcal{O}(\mathcal{F}')<\dim\mathcal{O}(\mathcal{F}).$
\end{rem}

\begin{rem}\label{nbptsing}
Si $\mathcal{F}$ d\'eg\'en\`ere sur $\mathcal{F}'$ le nombre de points
singuliers de $\mathcal{F}'$ (compt\'es sans multiplicit\'e) est plus 
petit ou \'egal \`a celui de $\mathcal{F}.$ En particulier si $\mathcal{F}$
n'a qu'un point singulier et d\'eg\'en\`ere sur $\mathcal{F}'$ alors
$\#\hspace{1mm}\mathrm{Sing}(\mathcal{F}')=1.$
\end{rem}

\begin{pro}
{\sl \textbf{\textit{1.}} Les orbites $\mathcal{O}(\mathcal{F}_1)$ et $\mathcal{O}(
\mathcal{F}_5)$ sont ferm\'ees.

\noindent\textbf{\textit{2.}} Pour $i=2$ et $3$ nous avons: $\overline{\mathcal{O}
(\mathcal{F}_i)}=\mathcal{O}(\mathcal{F}_i)\cup\mathcal{O}(\mathcal{F}_1).$}
\end{pro}

\begin{proof}[\sl D\'emonstration]
Le point \textbf{\textit{1.}} r\'esulte comme nous l'avons dit du fait que 
$\mathcal{O}(\mathcal{F}_1)$ et $\mathcal{O}(\mathcal{F}_5)$
r\'ealisent la dimension minimale.

\noindent La Remarque \ref{nbptsing} et des raisons de dimension 
font que les $\mathcal{F}_i,$ 
avec $i=2,$ $3,$ ne peuvent d\'eg\'en\'erer que sur~$\mathcal{F}_1.$
Montrons qu'il en est ainsi. Consid\'erons la famille d'automorphismes
$\varphi=\varphi_\varepsilon=\left(\frac{x}{\varepsilon^3},\frac{y}{
\varepsilon^2}\right).$ Nous avons
\begin{align*}
\varepsilon^9\varphi^*\omega_2=x^2\mathrm{d}x+(\varepsilon x+y^2)(x\mathrm{d}y-y\mathrm{d}x)
\end{align*}
qui tend vers $\omega_1$ lorsque $\varepsilon$ tend vers $0.$ Le
feuilletage $\mathcal{F}_2$ d\'eg\'en\`ere donc sur $\mathcal{F}_1.$

\noindent Pour montrer que $\mathcal{F}_3$ d\'eg\'en\`ere sur 
$\mathcal{F}_1$ nous proc\'edons comme suit. Nous 
consid\'erons la famille~$(\omega_\varepsilon)$ donn\'ee par 
$$\omega_\varepsilon=x(\varepsilon y+(1-\varepsilon)x)\mathrm{d}x+
(\varepsilon x^2+y^2)(x\mathrm{d}y-y\mathrm{d}x)$$
et d\'efinissant la famille de feuilletages $(\mathcal{G}_\varepsilon).$ Nous constatons
que
\begin{align*}
& \mathcal{G}_0=\mathcal{F}_1, &&\mathcal{G}_1=\mathcal{F}_3 &&
\text{et} &&\mu(\mathcal{G}_\varepsilon,0)=7.
\end{align*}
Par suite pour chaque $\varepsilon,$ \`a conjugaison pr\`es,
$\mathcal{G}_\varepsilon$ est l'un des $\mathcal{F}_i.$
Comme le $1$-jet de $\omega_\varepsilon$ est nul,~$\mathcal{G}_\varepsilon$
est conjugu\'e \`a $\mathcal{F}_1,$ $\mathcal{F}_2$ ou $\mathcal{F}_3.$
Pour des raisons de dimension $\mathcal{O}(\mathcal{F}_1)$ 
(resp. $\mathcal{O}(\mathcal{F}_2)$) ne peut adh\'erer \`a 
$\mathcal{O}(\mathcal{F}_3).$ Il en r\'esulte que pour $\varepsilon$ voisin
de $1$ le feuilletage $\mathcal{G}_\varepsilon$ appartient \`a
$\mathcal{O}(\mathcal{F}_3).$ Ainsi $\mathcal{G}_\varepsilon$
est dans $\mathcal{O}(\mathcal{F}_3)$ pour presque 
tout $\varepsilon$ et~$\mathcal{F}_3$ 
d\'eg\'en\`ere sur $\mathcal{F}_1.$ 
\end{proof}

\bigskip

\noindent Rappelons une notion introduite par \textsc{Pereira}~(\cite{Pe}).

\begin{defi}
Soient $\mathcal{F}$ un feuilletage sur $\mathbb{P}^2
(\mathbb{C})$ et $m$ un point r\'egulier de $\mathcal{F}.$ On
dit que $m$ est {\sl d'inflexion ordinaire} pour $\mathcal{F}$ si la
feuille $\mathcal{L}_m$ de $\mathcal{F}$ passant par $m$ a un 
point d'inflexion ordinaire en $m;$ notons $\mathrm{Flex}(\mathcal{F})$
l'adh\'erence de l'ensemble de ces points.
\end{defi}

\begin{egs}\label{egflex}
\begin{itemize}
\item[\texttt{1. }] Dans la carte $X=1$ le feuilletage $\mathcal{F}_1$ a pour
int\'egrale premi\`ere $\frac{y^3}{3}-z.$ En r\'esulte que 
$\mathrm{Flex}(\mathcal{F}_1)$ est la droite $y=0.$

\item[\texttt{2. }] Le feuilletage $\mathcal{F}_5$ a pour int\'egrale premi\`ere 
$\frac{y}{x}-\frac{1}{y};$ de sorte que ses feuilles ont pour 
adh\'erence des coniques, aussi $\mathrm{Flex}(\mathcal{F}_5)$
est vide. 
\end{itemize}
\end{egs} 

\bigskip

\noindent Soit $\mathcal{Z}=E\frac{\partial}{\partial X}+F\frac{\partial}
{\partial Y}+G\frac{\partial}{\partial Z}$ un champ de vecteurs homog\`ene
de degr\'e $2$ sur $\mathbb{C}^3$ non colin\'eaire au champ radial, 
$\text{pgcd}(E,F,G)=1;$ \`a 
un tel champ est associ\'e un feuilletage quadratique
$\mathcal{F}$ sur~$\mathbb{P}^2(\mathbb{C})$ d\'efini par la~$1$-forme
$\omega=i_\mathcal{R}i_\mathcal{Z}\mathrm{d}X\wedge \mathrm{d}Y\wedge \mathrm{d}Z.$ Introduisons
le polyn\^ome $\mathcal{H}$ de degr\'e $6$ 
$$\mathcal{H}(X,Y,Z)=\left\vert\begin{array}{ccc}
X & E & \mathcal{Z}(E)\\
Y & F & \mathcal{Z}(F)\\
Z & G & \mathcal{Z}(G)
\end{array}
\right\vert;$$
notons que $\mathcal{H}$ ne d\'epend pas du choix du champ
de vecteurs $\mathcal{Z}$ dans le noyau de $\omega.$
D'apr\`es \cite{Pe} le lieu des z\'eros de $\mathcal{H}$
est constitu\'e de $\mathrm{Flex}(\mathcal{F})$ et de 
l'ensemble des droites invariantes par $\mathcal{F}.$ 
\bigskip

\noindent L'\'enonc\'e qui suit montre qu'un feuilletage
quadratique g\'en\'eral $\mathcal{F}$ 
d\'eg\'en\`ere sur $\mathcal{F}_1.$

\begin{thm}\label{degenere}
{\sl Il existe un ensemble alg\'ebrique $\Sigma$ non trivial, contenu dans 
$\mathscr{F}(2;2)$ ayant la propri\'et\'e suivante: tout feuilletage de 
$\mathscr{F}(2;2)\setminus\Sigma$ d\'eg\'en\`ere sur $\mathcal{F}_1.$ 

\noindent En particulier, pour tout $\mathcal{F}$ dans $\mathscr{F}(2;2)
\setminus\Sigma,$ l'orbite $\mathcal{O}(\mathcal{F})$ de $\mathcal{F}$ n'est pas
ferm\'ee.}
\end{thm}

\noindent Donnons une interpr\'etation g\'eom\'etrique de l'\'enonc\'e qui
appara\^itra en cours de d\'emonstration. 
Un feuilletage~$\mathcal{F}$ appartient \`a $\mathscr{F}(2;2)\setminus~\Sigma$ 
si et seulement si l'ensemble $\mathrm{Flex}(\mathcal{F})$
est non vide. Cette propri\'et\'e est g\'en\'erique dans l'ensemble 
des feuilletages.

\begin{proof}[\sl D\'emonstration]
Consid\'erons la famille $\varphi_\varepsilon=\varphi=(\varepsilon^3x,
\varepsilon y).$ Notons que
\begin{itemize}
\item[\texttt{1. }] $\varphi^*(x^iy^j\mathrm{d}x)=\varepsilon^{3i+j+3}x^iy^j\mathrm{d}x$ est divisible 
par $\varepsilon^3$ et $\frac{1}{\varepsilon^3}\varphi^*(x^iy^j\mathrm{d}x)$ 
tend vers $0$ lorsque $\varepsilon$ tend vers $0$ sauf pour $i=j=0;$

\item[\texttt{2. }] $\varphi^*(x^iy^j\mathrm{d}y)=\varepsilon^{3i+j+1}x^iy^j\mathrm{d}y$ est divisible 
par $\varepsilon^3$ sauf pour $i=j=0$ et $(i,j)=(0,1).$ Si $(i,j)$ 
n'appartient pas \`a $\{(0,0),\hspace{1mm}(0,1),\hspace{1mm}(0,2)\},$ la forme
$\frac{1}{\varepsilon^3}\varphi^*(x^iy^j\mathrm{d}y)$ tend vers $0$ lorsque
$\varepsilon$ tend vers $0.$
\end{itemize}
Si $\mathcal{F}$ est d\'efini par une $1$-forme $\omega$
du type
\begin{align*}
&(\alpha +*x+*y+*x^2+*xy+*y^2)\mathrm{d}x+(*x+*x^2+*xy+\beta y^2)\mathrm{d}y\\
&\hspace{3mm} +(*x^2+*xy+*y^2)(y\mathrm{d}x-x\mathrm{d}y), &&\text{o\`u les }*\in\mathbb{C},\hspace{1mm}
\alpha,\hspace{1mm}\beta\in\mathbb{C}^*,
\end{align*}
nous dirons que $\omega$ ou $\mathcal{F}$ v\'erifie la propri\'et\'e $(\mathcal{P}).$ Si tel est 
le cas nous constatons que
$$\lim_{\varepsilon\to 0}\frac{1}{\varepsilon^3}\varphi^*\omega=\alpha\mathrm{d}x+
\beta y^2\mathrm{d}y.$$
Visiblement $\alpha\mathrm{d}x+\beta y^2\mathrm{d}y$ d\'efinit un feuilletage 
conjugu\'e \`a $\mathcal{F}_1;$
en particulier si $\mathcal{F}$ v\'erifie $(\mathcal{P})$ alors $\mathcal{F}$ d\'eg\'en\`ere sur $\mathcal{F}_1.$

\noindent Notons que si $\mathcal{F}$ a la propri\'et\'e $(\mathcal{P})$ le point $(0,0)$ est 
non singulier (car~$\alpha\not=0$) et la feuille $\mathcal{L}_{(0,0)}$
de $\mathcal{F}$ en $(0,0)$ est tangente \`a l'axe $x=0.$ Comme le coefficient
de $\omega$ sur $y\mathrm{d}y$ est nul un calcul \'el\'ementaire montre que 
$\mathcal{L}_{(0,0)}$ a un point d'inflexion en~$(0,0).$ Finalement le fait
que $\beta$ soit non nul implique que ce point d'inflexion est d'ordre minimal $1.$

\noindent Soit $\mathcal{F}$ un feuilletage quadratique sur 
$\mathbb{P}^2(\mathbb{C});$ supposons que $\mathrm{Flex}(\mathcal{F})$
soit non vide. Choisissons alors un 
syst\`eme de coordonn\'ees $(x,y)$ tel que $m=(0,0)$ soit
d'inflexion ordinaire et $x=0$ soit la tangente \`a~$\mathcal{L}_m$ en~$m.$ La $1$-forme $\omega$
d\'ecrivant $\mathcal{F}$ satisfait la propri\'et\'e $(\mathcal{P})$ et $\mathcal{F}$ d\'eg\'en\`ere sur $\mathcal{F}_1.$

\noindent L'ensemble $\Sigma$ est d\'efini comme l'ensemble des feuilletages $\mathcal{F}$ tel que $\mathrm{Flex}(\mathcal{F})=\emptyset.$
\end{proof}

\begin{rems}
\begin{itemize}
\item[\texttt{1. }] Comme indiqu\'e dans l'exemple \ref{egflex} le feuilletage 
$\mathcal{F}_1$ poss\`ede une droite de points d'inflexion qui sont d'ordre $1;$ ceci 
confirme que $\Sigma$ est propre.

\item[\texttt{2. }]
L'ensemble $\Sigma$ est invariant sous l'action de $\mathrm{PGL}_3(
\mathbb{C}).$

\item[\texttt{3. }] De la m\^eme fa\c{c}on que pour $\mathcal{F}_5,$ les 
feuilletages ayant une int\'egrale premi\`ere du type $Q_1/Q_2,$
avec $Q_1,$ $Q_2$ polyn\^ome de degr\'e inf\'erieur ou \'egal \`a $2,$
ne poss\`edent pas de points d'inflexion car les feuilles sont d'adh\'erence
des coniques; ils appartiennent \`a $\Sigma.$ 
\end{itemize}
\end{rems}

\noindent On peut voir ais\'ement que l'ensemble $\mathrm{Flex}
(\mathcal{F}_4)$ est non vide comme le montrent d'ailleurs les figures 
du \S\,\ref{portraitdephase}. 

\begin{cor}
{\sl Le feuilletage $\mathcal{F}_4$ d\'eg\'en\`ere sur $\mathcal{F}_1;$
en particulier $\mathcal{O}(\mathcal{F}_4)$ n'est pas ferm\'ee.}
\end{cor}

\begin{proof}[\sl D\'emonstration]
Le feuilletage $\mathcal{F}_4$ est d\'ecrit par 
$$\omega_4=(x+y^2-x^2y)\mathrm{d}y+x(x+y^2)\mathrm{d}x.$$
Soit $\alpha$ dans $\mathbb{C}$ tel que $\alpha^3=4;$ posons
$\beta=\frac{\alpha^2}{2}.$ En faisant agir la translation 
$(x+\alpha,y+\beta)$ sur~$\omega_4$ nous obtenons une $1$-forme
du type
\begin{small}
$$\omega=(\alpha(\alpha+\beta^2)+*x+*y+*x^2+*xy+*y^2)\mathrm{d}x+(*x+*x^2+*xy+y^2)\mathrm{d}y+
(*x^2+*xy+*y^2)(y\mathrm{d}x-x\mathrm{d}y);$$
\end{small}
la propri\'et\'e $(\mathcal{P}),$
introduite au Th\'eor\`eme \ref{degenere}, est donc satisfaite. 
D'apr\`es la d\'emonstration du Th\'eor\`eme \ref{degenere} nous avons
$$\omega_0=\lim_{\varepsilon\to 0}\frac{1}{\varepsilon^3}(\varepsilon^3x,\varepsilon y)^*\omega=
\alpha(\alpha+\beta^2) \mathrm{d}x+y^2\mathrm{d}y.$$

\noindent On v\'erifie que $\alpha(\alpha+\beta^2)$ est non nul ce 
qui implique que le feuilletage induit par $\omega_0$ est 
conjugu\'e \`a $\mathcal{F}_1.$
\end{proof}

\bigskip

\begin{eg}[de \textsc{Jouanolou}]
Introduisons le feuilletage $\mathcal{F}_J$\label{not8b} de degr\'e 
$2$ sur $\mathbb{P}^2(\mathbb{C})$ d\'efini, dans la carte affine~$Z=1,$ par
$$\omega_J=(x^2y-1)\mathrm{d}x+(y^2-x^3)\mathrm{d}y.$$\label{not8c}
Cet exemple, d\^u \`a \textsc{Jouanolou}, est tr\`es populaire. C'est
historiquement le premier exemple explicite 
de feuilletage sans courbe alg\'ebrique invariante (\cite{J}); c'est aussi
un feuilletage qui n'admet pas d'ensemble minimal non trivial (\cite{CdF}). 
La proc\'edure appliqu\'ee dans la d\'emonstration du Th\'eor\`eme
\ref{degenere} s'applique directement, autrement dit la propri\'et\'e 
$(\mathcal{P})$ est v\'erifi\'ee, ainsi $\mathcal{F}_J$ d\'eg\'en\`ere sur~$\mathcal{F}_1.$
\end{eg}

\medskip

\noindent Pour $\mathcal{F}$ g\'en\'erique dans $\mathscr{F}(2;2)$ 
on peut se demander si $\mathcal{F}$ d\'eg\'en\`ere uniquement 
sur $\mathcal{F}_1.$ Les consid\'erations suivantes vont nous 
montrer que non. 

\noindent Soit $\mathcal{F}$ un \'el\'ement de $\mathscr{F}(2;2)$ 
d\'ecrit par une $1$-forme $\omega.$ Supposons
que $(0,0)$ soit un point singulier de $\mathcal{F};$ nous pouvons
trouver des coordonn\'ees o\`u $\omega$ est du type
\begin{align*}
& (\alpha y+*x+*x^2+*y^2+*xy)\mathrm{d}x+(\beta x+\gamma y^2+*x^2+*xy)\mathrm{d}y &&\\
& \hspace{3mm}+(*x^2+*xy+*y^2)(y\mathrm{d}x-x\mathrm{d}y), && \text{avec } *,\, \alpha,\hspace{1mm}\beta,
\hspace{1mm}\gamma\in\mathbb{C}.
\end{align*}
Posons $\varphi=(\varepsilon^2x,\varepsilon y).$ Nous constatons
que
\begin{align*}
\Omega_0=\lim_{\varepsilon\to 0}\frac{1}{\varepsilon^3}\varphi^*\omega=\alpha 
y\mathrm{d}x+(\beta x+\gamma y^2)\mathrm{d}y.
\end{align*}

\noindent Lorsque $\alpha,$ $\beta$ et $\gamma$ sont tous les trois non nuls 
$\Omega_0$ d\'efinit un feuilletage quadratique $\mathscr{F}_0$
sur lequel d\'eg\'en\`ere \'evidemment $\mathcal{F}.$ Alors que
la non nullit\'e des coefficients $\alpha$ et $\beta$ signifie que $(0,0)$
est de multiplicit\'e $1,$ celle de $\gamma$ s'interpr\`ete comme suit:
si $\gamma$ est nul $\mathcal{F}$ poss\`ede une droite invariante, 
ce qui g\'en\'eriquement n'arrive pas dans $\mathscr{F}(2;2).$

\noindent D\'etaillons le cas $\alpha\beta\gamma\not=0.$ On peut supposer 
\`a conjugaison pr\`es que $\mathcal{F}_0=\mathcal{F}_0
(\lambda)$ est d\'efini par 
\begin{align*}
&\Omega_0=\Omega_0(\lambda)=y\mathrm{d}x+(\lambda x+y^2)\mathrm{d}y,
&&\lambda=\beta/\alpha.
\end{align*}
Notons que l'invariant $\lambda$ est reli\'e \`a la singularit\'e
$(0,0)$ du feuilletage initial $\mathcal{F};$ plus pr\'ecis\'ement
l'invariant de \textsc{Baum}-\textsc{Bott} $\mathrm{BB}(\mathcal{F};
(0,0))$\label{not11} de $\mathcal{F}$ au point $(0,0)$ est \'egal \`a celui
de $\mathcal{F}_0(\lambda)$
$$\mathrm{BB}(\mathcal{F};(0,0))=-\frac{(\alpha-\beta)^2}{\alpha\beta}=-\frac{(1-
\lambda)^2}{\lambda}=\mathrm{BB}(\mathcal{F}_0(\lambda);
(0,0)).$$
Nous allons d\'ecrire les feuilletages $\mathcal{F}_0(\lambda)$
pour $\lambda$ dans $\mathbb{C}^*.$ Ils pr\'eservent les 
droites $y=0$ et $z=0$ et comptent deux points singuliers
$(1:0:0)$ et $(0:0:1).$ Le point $(0:0:1)$ est de multiplicit\'e~$1$ 
et par suite $(1:0:0)$ de multiplicit\'e $6.$ Nous constatons
que pour $\lambda=-1$ la~$1$-forme
$$\frac{\Omega_0(-1)}{y^2}=\frac{y\mathrm{d}x-x\mathrm{d}y}{y^2}+\mathrm{d}y$$
est ferm\'ee et a pour primitive $\frac{x}{y}+y.$ Il s'en suit que
$\mathcal{F}_0(-1)$ est conjugu\'e au feuilletage $\mathcal{F}_5.$
Il y a encore une valeur sp\'eciale de $\lambda$ qui intervient:
$\lambda=-2.$ Dans ce cas nous remarquons que~$\frac{\Omega_0
(-2)}{y^3}$ est ferm\'ee et s'int\`egre en $\frac{x}{y^2}+\ln y;$ ceci
permet de constater que $\mathrm{Iso}(\mathcal{F}_0(-2))$
est le groupe~$\{(\varepsilon^2x,\varepsilon y)\hspace{1mm}
\vert\hspace{1mm}\varepsilon\in\mathbb{C}^*\}.$ En particulier~$\mathcal{O}(\mathcal{F}_0(-2))$ est de dimension $7.$

\noindent Dans la suite nous supposons que $\lambda$
appartient \`a $\mathbb{C}\setminus\{0,\hspace{1mm} -1,
\hspace{1mm} -2\}.$ Un calcul \'el\'ementaire montre que la conique 
$(2+\lambda)x+y^2=0$ est invariante par $\mathcal{F}_0(
\lambda).$ Mieux la $1$-forme rationnelle~$\frac{\Omega_0(\lambda)}
{y((2+\lambda)x+y^2)}$ est ferm\'ee. Un calcul \'el\'ementaire en 
donne une primitive $$\ln((2+\lambda)x+y^2)+\lambda\ln y,$$
ce qui permet de v\'erifier que le groupe d'isotropie est ici encore
$$\mathrm{Iso}(\mathcal{F}_0(\lambda))=\{(\varepsilon^2x,\varepsilon
y)\hspace{1mm}\vert\hspace{1mm}\varepsilon\in\mathbb{C}^*\}.$$

\noindent Remarquons que deux feuilletages $\mathcal{F}_0(\lambda)$
et $\mathcal{F}_0(\lambda')$ sont conjugu\'es si et seulement si 
$\lambda=\lambda'.$

\noindent Des consid\'erations qui pr\'ec\`edent nous tirons l'\'enonc\'e
suivant.

\begin{thm}\label{sansdteinv}
{\sl Soit $\mathcal{F}$ un \'el\'ement de $\mathscr{F}(2;2)$ sans droite invariante.

\noindent\textbf{\textit{1.}} Le feuilletage $\mathcal{F}$ d\'eg\'en\`ere sur
$\mathcal{F}_5$ si et seulement si $\mathcal{F}$ poss\`ede un 
point singulier $m$ de multiplicit\'e $1$ v\'erifiant $\mathrm{BB}
(\mathcal{F};m)=4.$

\noindent\textbf{\textit{2.}} Si $\mathcal{F}$ poss\`ede un point 
singulier $m$ de multiplicit\'e $1$ tel que $\mathrm{BB}(\mathcal{F};
m)=2-\lambda-\frac{1}{\lambda},$ alors $\mathcal{F}$ 
d\'eg\'en\`ere sur~$\mathcal{F}_0(\lambda)$ donn\'e par
$$y\mathrm{d}x+(\lambda x+y^2)\mathrm{d}y.$$
En particulier si 
$\mathcal{F}$ est g\'en\'erique, $\mathcal{F}$ d\'eg\'en\`ere au
moins sur sept feuilletages (non conjugu\'es) \`a orbites de
dimension $7.$ }
\end{thm}

\begin{proof}[\sl D\'emonstration]
Concernant la premi\`ere assertion, il nous reste seulement \`a
prouver que si $\mathcal{F}$ d\'eg\'en\`ere sur~$\mathcal{F}_5,$
alors $\mathcal{F}$ poss\`ede un point singulier de multiplicit\'e
$1$ dont l'invariant de \textsc{Baum}-\textsc{Bott} est~$4.$ Si 
tel est le cas il y a une famille analytique $(\mathcal{F}_\varepsilon)$ 
d\'efinie par la famille de $1$-formes $(\omega_\varepsilon)$ 
telle que pour $\varepsilon$ non nul $\mathcal{F}_\varepsilon$
soit dans $\mathcal{O}(\mathcal{F})$ et pour $\varepsilon$ nul
nous ayons $\mathcal{F}_{\varepsilon=0}=\mathcal{F}_5.$ Le point singulier
de multiplicit\'e~$1,$ not\'e~$m_0,$ de $\mathcal{F}_5$ est 
\og stable\fg; il existe une famille analytique 
$(m_\varepsilon)$ de points singuliers de $\mathcal{F}_\varepsilon$ de
multiplicit\'e $1$ telle que $m_{\varepsilon=0}=~m_0.$ Les 
$\mathcal{F}_\varepsilon$ \'etant conjugu\'es pour $\varepsilon$ non nul,
$\mathrm{BB}(\mathcal{F}_\varepsilon;m_\varepsilon)$ est 
localement constant; par sui\-te~$\mathrm{BB}(\mathcal{F}_\varepsilon;
m_\varepsilon)=4$ pour $\varepsilon$ petit. En particulier 
$\mathcal{F}$ poss\`ede un point singulier comme dans 
l'\'enonc\'e.\medskip

\noindent Le point \textbf{\textit{2.}} est cons\'equence du fait 
suivant: g\'en\'eriquement $\mathcal{F}$ poss\`ede sept points
singuliers distincts \`a invariants de \textsc{Baum}-\textsc{Bott}
distincts.
\end{proof}

\begin{rem}
Les Th\'eor\`emes \ref{degenere} et \ref{sansdteinv} assurent en particulier
qu'un feuilletage g\'en\'erique d\'eg\'en\`ere sur 
$\mathcal{F}_1$ mais pas sur $\mathcal{F}_5.$
\end{rem}

\noindent Le Th\'eor\`eme \ref{sansdteinv} semble 
occulter les \'elements $\mathcal{F}$ de $\mathscr{F}(2;2)$ avec 
droite invariante. En fait si $\mathcal{F}$ poss\`ede 
un point singulier simple (condition $\alpha\beta\not=0$)
par lequel ne passe pas de droite invariante la
proc\'edure de d\'eg\'en\'erescence sur un 
$\mathcal{F}_0(\lambda)$ s'applique et autorise la pr\'esence de droites invariantes
ne passant pas par $(0,0).$ Une autre 
fa\c{c}on de proc\'eder est la suivante. Si $\mathcal{F}$
laisse $y=0$ invariante la $1$-forme qui le d\'efinit s'\'ecrit 
$$\omega=ya(x,y)\mathrm{d}x+b(x,y)\mathrm{d}y.$$
En faisant agir la transformation lin\'eaire diagonale $(x,\varepsilon y)$ sur $\omega$ 
puis en passant \`a la limite lorsque $\varepsilon\to 0$ on obtient
$$\omega_0=ya(x,0)\mathrm{d}x+b(x,0)\mathrm{d}y$$
qui g\'en\'eriquement sur $\mathcal{F}$ d\'efinit un feuilletage
de degr\'e $2$ sur lequel $\mathcal{F}$ d\'eg\'en\`ere.

\begin{pro}
{\sl L'ensemble $\mathscr{F}(2;2)$ contient 
des \'el\'ements qui d\'eg\'en\`erent \`a la fois sur 
$\mathcal{F}_1$ et~$\mathcal{F}_5.$}
\end{pro}

\begin{proof}[\sl D\'emonstration]
Rappelons que la $1$-forme 
$$\omega_1=x^2\mathrm{d}x+y^2(x\mathrm{d}y-y\mathrm{d}x)$$
d\'efinit $\mathcal{F}_1,$ feuilletage pour lequel l'ensemble
$\mathrm{Flex}(\mathcal{F}_1)$ est une droite. Un feuilletage $\mathcal{F}_\varepsilon$
proche de $\mathcal{F}_1$ va donc satisfaire $\mathrm{Flex}
(\mathcal{F}_\varepsilon)\not=\emptyset$ et par suite d\'eg\'en\`ere
sur $\mathcal{F}_1$ (Th\'eor\`eme \ref{degenere}). En particulier il en est 
ainsi pour $\mathcal{F}_\varepsilon$
d\'efini par
$$\omega_\varepsilon=\varepsilon(x\mathrm{d}y-y\mathrm{d}x)+\omega_1$$
tout du moins pour $\varepsilon$ petit. Mais pour $\varepsilon$
non nul le $1$-jet de $\omega_\varepsilon$ en $(0,0)$ est 
$x\mathrm{d}y-y\mathrm{d}x$ et le coefficient de $x^2\mathrm{d}x$ est $1$ ce qui nous 
permet de faire d\'eg\'en\'erer (\`a conjugaison pr\`es)
$\mathcal{F}_\varepsilon$ sur $\mathcal{F}_5$ pour 
$\varepsilon\not=0;$ ainsi pour $\varepsilon$ petit non nul
$\mathcal{F}_\varepsilon$ d\'eg\'en\`ere \`a la fois sur
$\mathcal{F}_1$ et~$\mathcal{F}_5.$
\end{proof}

\noindent Soit $\mathcal{F}$ dans $\mathscr{F}(2;2);$ nous
avons l'alternative
\begin{itemize}
\item[\texttt{1. }] si $\mathrm{Flex}(\mathcal{F})$ est non vide, 
$\mathcal{F}$ d\'eg\'en\`ere sur $\mathcal{F}_1$
(Th\'eor\`eme \ref{degenere});

\item[\texttt{2. }] si $\mathrm{Flex}(\mathcal{F})=\emptyset$ alors 
$\mathcal{F}$ ne d\'eg\'en\`ere pas sur $\mathcal{F}_1$
(stabilit\'e des points d'inflexion).
\end{itemize}

\smallskip

\noindent Parmi les \'el\'ements
de $\mathscr{F}(2;2)$ sans point d'inflexion, il y a les
feuilletages d\'efinis par un pinceau de coniques. 
Leur classification est connue
depuis longtemps; on peut par exemple la trouver dans 
\cite{HP}, chapitre XIII, \S 11. Elle correspond \`a celle 
des feuilletages de degr\'e inf\'erieur ou
\'egal \`a $2$ ayant une int\'egrale premi\`ere rationnelle
du type $Q_1/Q_2,$ les $Q_i$ d\'esignant des formes
quadratiques \`a $3$ variables, le pinceau associ\'e
\'etant $Q_1+\lambda Q_2.$ Il y a cinq mod\`eles \`a conjugaison pr\`es;
deux conduisent \`a des feuilletages de degr\'e~$1$
(ils sont d\'ecrits par $xy=$ cte, $y^2-x=$ cte), les 
trois autres sont les suivants: le feuilletage~$\mathcal{F}_6$\label{not14}
qui repr\'esente le pinceau g\'en\'erique et qui a pour
int\'egrale premi\`ere $\frac{x(y+z)}{y(x+z)},$ le feuilletage~$\mathcal{F}_5$ d\'ej\`a rencontr\'e et enfin le feuilletage
$\mathcal{F}_7$\label{not15} dont voici une int\'egrale premi\`ere
$\frac{xz}{y(y-x)}.$
L'ensemble des coniques dans~$\mathbb{P}^2(\mathbb{C})$
est isomorphe \`a $\mathbb{P}^5(\mathbb{C})$ (\emph{voir} 
\cite{H}); il en r\'esulte qu'un pinceau de coniques s'identifie
\`a une droite dans $\mathbb{P}^5(\mathbb{C}),$ {\it i.e.}
\`a un point dans la grassmanienne $\mathbb{G}(1;5)$ 
des droites de $\mathbb{P}^5(\mathbb{C}),$ qui est de 
dimension $8.$ Ainsi~$\mathrm{PGL}_3(\mathbb{C})$
agit sur $\mathbb{G}(1;5)$ et l'orbite du pinceau 
g\'en\'erique est ouverte et dense. Par suite $\mathcal{O}
(\mathcal{F}_6)$ adh\`ere \`a toutes les autres orbites 
(associ\'ees \`a des pinceaux), en particulier \`a celle de 
$\mathcal{F}_5$ et \`a celle de $\mathcal{F}_7.$
Une mani\`ere, autre que la mani\`ere classique, de 
distinguer ces trois feuilletages est de consid\'erer 
leur groupe d'isotropie: celui de $\mathcal{F}_5$ est, 
comme nous l'avons vu, 
$$\mathrm{Iso}(\mathcal{F}_5)=\left\{(\alpha^2 x,\alpha y),
\left(\frac{x}{1+\beta y},\frac{y}{1+\beta y}\right)
\hspace{1mm}\Big\vert\hspace{1mm}\alpha\in\mathbb{C}^*,
\hspace{1mm}\beta\in\mathbb{C}\right\},$$
celui de $\mathcal{F}_6$ est fini et celui de $\mathcal{F}_7$
est donn\'e par
$$\mathrm{Iso}(\mathcal{F}_7)=\{(\alpha x,\alpha y)
\hspace{1mm}\vert\hspace{1mm}\alpha\in\mathbb{C}^*\}.$$

\noindent Tout ceci est r\'esum\'e dans la:

\begin{pro}
{\sl Soit $\mathcal{F}$ un \'el\'ement de $\mathscr{F}(2,2)$ associ\'e \`a un pinceau de 
coniques. Alors $\mathcal{F}$ est conjugu\'e \`a $\mathcal{F}_5,$ ou $\mathcal{F}_6$
ou $\mathcal{F}_7.$ De plus, les adh\'erences \'etant prises dans $\mathscr{F}(2;2),$
nous avons 
\smallskip
\begin{itemize}
\item[\texttt{1. }] $\overline{\mathcal{O}(\mathcal{F}_5)}= \mathcal{O}(\mathcal{F}_5)$
et $\dim\mathcal{O}(\mathcal{F}_5)=6;$

\item[\texttt{2. }] $\overline{\mathcal{O}(\mathcal{F}_7)}=
\mathcal{O}(\mathcal{F}_5)\cup\mathcal{O}(\mathcal{F}_7)$ et 
$\dim\mathcal{O}(\mathcal{F}_7)=7;$

\item[\texttt{3. }] $\overline{\mathcal{O}(\mathcal{F}_6)}=\mathcal{O}(\mathcal{F}_5)\cup
\mathcal{O}(\mathcal{F}_6)\cup\mathcal{O}(\mathcal{F}_7)$ et 
$\dim\mathcal{O}(\mathcal{F}_6)=8.$ 
\end{itemize}}
\end{pro}

\begin{proof}[\sl D\'emonstration]
Une fois que l'on a constat\'e que si $\mathcal{F}$ d\'eg\'en\`ere
sur $\mathtt{F}$ alors $\mathtt{F}$ a aussi une int\'egrale
premi\`ere du type $Q_1/Q_2$ il nous reste seulement \`a d\'emontrer
que $\mathcal{F}_7$ d\'eg\'en\`ere sur $\mathcal{F}_5.$
Introduisons le feuilletage $\mathcal{F}'$ donn\'e par les niveaux de 
$$\frac{y}{x}+\frac{x}{y}+\frac{1}{y};$$
son groupe d'isotropie est 
$$\mathrm{Iso}(\mathcal{F}')=\left\{\left(\frac{x}{1+\alpha y},
\frac{y}{1+\alpha y}\right)\big\vert\alpha\in\mathbb{C}^*\right\}.$$
En particulier $\dim\mathcal{O}(\mathcal{F}')=7$ et $\mathcal{F}'$ 
appartient \`a $\mathcal{O}(\mathcal{F}_7).$
Quitte \`a faire agir $\varphi_\varepsilon=(\varepsilon^2x,
\varepsilon y)$ sur~$\frac{y}{x}~+~\frac{x}{y}+~\frac{1}{y}$ nous
obtenons \`a multiplication par une constante pr\`es
$$\varepsilon^2\frac{x}{y}+\frac{y}{x}+\frac{1}{y};$$
nous retrouvons lorsque $\varepsilon$ tend vers $0$ 
l'int\'egrale premi\`ere de $\mathcal{F}_5.$
\end{proof}

\noindent Tout ceci laisse penser que les seules orbites 
ferm\'ees dans $\mathscr{F}(2;2)$ sous l'action de 
$\mathrm{Aut}(\mathbb{P}^2(\mathbb{C}))$ sont celles
de $\mathcal{F}_1$ et $\mathcal{F}_5.$

\subsection{En dimension sup\'erieure}\hspace{1mm}

\noindent L'\'etude de la nature des orbites sous l'action du groupe 
lin\'eaire pour les feuilletages de dimension sup\'erieure \`a $2$ se
pose aussi naturellement, tout comme la classification des feuilletages
ayant des \og petites\fg\hspace{1mm} singularit\'es.

\noindent Soient $\mathcal{F}_0$ dans $\mathscr{F}(2;N)$
et $\mathcal{F}$ dans $\mathscr{F}(n;N)$ (o\`u $\mathscr{F}(n;N)$ d\'esigne l'ensemble des feuilletages
de degr\'e $N$ sur l'espace projectif~$\mathbb{P}^n(\mathbb{C})$). On dit que $\mathcal{F}$ est un 
{\sl d\'eploiement} de $\mathcal{F}_0$ s'il existe un plonge\-ment~$i\hspace{1mm}
\colon\hspace{1mm}\mathbb{P}^2(\mathbb{C})\to~\mathbb{P}^n(\mathbb{C})$
lin\'eaire tel que~$i^*\mathcal{F}=\mathcal{F}_0.$

\noindent Rappelons (\emph{voir} \cite{CLN2}) que si $\mathcal{F}$ est un feuilletage quadratique 
sur $\mathbb{P}^n(\mathbb{C})$ alors
\begin{itemize}
\item[\texttt{1. }] ou bien $\mathcal{F}$ est d\'efini par une $1$-forme ferm\'ee;

\item[\texttt{2. }] ou bien $\mathcal{F}=\pi^*\mathcal{F}_0$ o\`u $\mathcal{F}_0$
d\'esigne un feuilletage quadratique de $\mathbb{P}^2(\mathbb{C})$ et
$\pi$ la projection de~$\mathbb{P}^n(\mathbb{C})$ sur~$\mathbb{P}^2(\mathbb{C}).$ 
\end{itemize}
\smallskip

\noindent Nous en d\'eduisons la:

\begin{pro}\label{depl}
{\sl Soit $\mathcal{F}_0$ un feuilletage de degr\'e~$2$ sur 
$\mathbb{P}^2(\mathbb{C}).$
Si $\mathcal{F}_0$ n'est pas d\'efini par une forme ferm\'ee tout d\'eploiement
$\mathcal{F}$ de $\mathcal{F}_0$ est trivial, {\it i.e.} il existe $\pi\hspace{1mm}
\colon\hspace{1mm}\mathbb{P}^n(\mathbb{C})\to \mathbb{P}^2(\mathbb{C})$
tel que~$\mathcal{F}=\pi^*\mathcal{F}_0.$}
\end{pro}

\noindent Le feuilletage $\mathcal{F}_4$ satisfait la Proposition \ref{depl}; 
ainsi $\mathcal{F}_4$ va jouer un r\^ole particulier dans l'\'etude des 
orbites des feuilletages quadratiques de $\mathbb{P}^n(\mathbb{C})$ avec 
$n\geq 3$ sous l'action de $\mathrm{Aut}(\mathbb{P}^n(\mathbb{C})).$

\section{Quelques portraits de phase en r\'eel}\label{portraitdephase}\hspace{1mm}

\noindent Les feuilletages $\mathcal{F}_i$ et $\mathcal{F}_J$ sont d\'efinis par des formes
\`a coefficients r\'eels (m\^eme entiers) et induisent donc des feuilletages 
de $\mathbb{R}^2$ dont nous proposons le portrait de phase; ces
portraits de phase ont \'et\'e trac\'es \`a l'aide d'un logiciel disponible sur 
http://www.math.psu.edu/melvin/phase/newphase.html.

\bigskip

\begin{center}
\begin{tabular}{cc}
\begin{picture}(0,0)%
\includegraphics{omega1_1.pstex}%
\end{picture}%
\setlength{\unitlength}{3947sp}%
\begingroup\makeatletter\ifx\SetFigFont\undefined%
\gdef\SetFigFont#1#2#3#4#5{%
  \reset@font\fontsize{#1}{#2pt}%
  \fontfamily{#3}\fontseries{#4}\fontshape{#5}%
  \selectfont}%
\fi\endgroup%
\begin{picture}(2400,2400)(1201,-7561)
\end{picture}%
\hspace{2cm}&\hspace{2cm}\begin{picture}(0,0)%
\includegraphics{omega2_1.pstex}%
\end{picture}%
\setlength{\unitlength}{3947sp}%
\begingroup\makeatletter\ifx\SetFigFont\undefined%
\gdef\SetFigFont#1#2#3#4#5{%
  \reset@font\fontsize{#1}{#2pt}%
  \fontfamily{#3}\fontseries{#4}\fontshape{#5}%
  \selectfont}%
\fi\endgroup%
\begin{picture}(2400,2400)(1201,-7561)
\end{picture}%
\\
& \\
$\mathcal{F}_1$ & \hspace{2cm} $\mathcal{F}_2$
\end{tabular}
\end{center}

\bigskip

\begin{center}
\begin{tabular}{cc}
\begin{picture}(0,0)%
\includegraphics{omega3_3.pstex}%
\end{picture}%
\setlength{\unitlength}{3947sp}%
\begingroup\makeatletter\ifx\SetFigFont\undefined%
\gdef\SetFigFont#1#2#3#4#5{%
  \reset@font\fontsize{#1}{#2pt}%
  \fontfamily{#3}\fontseries{#4}\fontshape{#5}%
  \selectfont}%
\fi\endgroup%
\begin{picture}(2400,2400)(1201,-7561)
\end{picture}%
\hspace{2cm}&\hspace{2cm}\begin{picture}(0,0)%
\includegraphics{jouanolou.pstex}%
\end{picture}%
\setlength{\unitlength}{3947sp}%
\begingroup\makeatletter\ifx\SetFigFont\undefined%
\gdef\SetFigFont#1#2#3#4#5{%
  \reset@font\fontsize{#1}{#2pt}%
  \fontfamily{#3}\fontseries{#4}\fontshape{#5}%
  \selectfont}%
\fi\endgroup%
\begin{picture}(2400,2400)(1201,-7561)
\end{picture}%
\\
& \\
$\mathcal{F}_3$ & \hspace{2cm}$\mathcal{F}_J$
\end{tabular}
\end{center}

\bigskip

\noindent Soit $\widetilde{\mathcal{F}}$ le feuilletage d\'ecrit 
dans la carte affine $Z=1$ par la 
$1$-forme $y^7\mathrm{d}x-x\mathrm{d}y.$ Dans le contexte r\'eel les feuilletages 
$\widetilde{\mathcal{F}}$ et $\mathcal{F}_4$ sont localement 
diff\'erentiablement conjugu\'es en l'origine. Ainsi le portrait de 
phase de~$\widetilde{\mathcal{F}}$ explique celui de $\mathcal{F}_4$
en $0;$ on observe un ph\'enom\`ene de type \og fleuve\fg\hspace{1mm}
\`a l'origine (\cite{D}).

\bigskip

\begin{center}
\begin{tabular}{cc}
\begin{picture}(0,0)%
\includegraphics{omega4_3.pstex}%
\end{picture}%
\setlength{\unitlength}{3947sp}%
\begingroup\makeatletter\ifx\SetFigFont\undefined%
\gdef\SetFigFont#1#2#3#4#5{%
  \reset@font\fontsize{#1}{#2pt}%
  \fontfamily{#3}\fontseries{#4}\fontshape{#5}%
  \selectfont}%
\fi\endgroup%
\begin{picture}(2400,2400)(1201,-7561)
\end{picture}%
\hspace{2cm}&\hspace{2cm}\begin{picture}(0,0)%
\includegraphics{modele.pstex}%
\end{picture}%
\setlength{\unitlength}{3947sp}%
\begingroup\makeatletter\ifx\SetFigFont\undefined%
\gdef\SetFigFont#1#2#3#4#5{%
  \reset@font\fontsize{#1}{#2pt}%
  \fontfamily{#3}\fontseries{#4}\fontshape{#5}%
  \selectfont}%
\fi\endgroup%
\begin{picture}(2400,2400)(1201,-7561)
\end{picture}%
\\
& \\
$\mathcal{F}_4$ & \hspace{2cm}$\widetilde{\mathcal{F}}$
\end{tabular}
\end{center}

\noindent Le portrait de phase de $\mathcal{F}_4$ ci-dessus 
sugg\`ere que le feuilletage r\'eel induit par $\mathcal{F}_4$
ne poss\`ede pas de cycle limite, plus pr\'ecis\'ement que toute
trajectoire adh\`ere au point singulier.

\bigskip

\noindent Pour avoir une meilleure id\'ee des trajectoires r\'eelles du 
feuilletage $\mathcal{F}_4$ consid\'erons le feuilletage~$\mathcal{F}_4^\perp$\label{not12}
donn\'e dans la carte affine $Z=1$ par la $1$-forme
$$\omega_4^\perp=-(x+y^2-x^2y)\mathrm{d}x+x(x+y^2)\mathrm{d}y.$$
\label{not13}

\noindent C'est un feuilletage de degr\'e $3$ dont les trajectoires 
sont dans la carte $Z=1$ orthogonales \`a celles de $\mathcal{F}_4.$
Nous proposons des repr\'esentations de $\mathcal{F}_4$ et 
$\mathcal{F}_4^\perp$ dans des domaines de tailles diff\'erentes

\bigskip

\begin{center}
\begin{tabular}{cc}
\begin{picture}(0,0)%
\includegraphics{omega4_3.pstex}%
\end{picture}%
\setlength{\unitlength}{3947sp}%
\begingroup\makeatletter\ifx\SetFigFont\undefined%
\gdef\SetFigFont#1#2#3#4#5{%
  \reset@font\fontsize{#1}{#2pt}%
  \fontfamily{#3}\fontseries{#4}\fontshape{#5}%
  \selectfont}%
\fi\endgroup%
\begin{picture}(2400,2400)(1201,-7561)
\end{picture}%
\hspace{2cm}&\hspace{2cm}\begin{picture}(0,0)%
\includegraphics{orthogonal1.pstex}%
\end{picture}%
\setlength{\unitlength}{3947sp}%
\begingroup\makeatletter\ifx\SetFigFont\undefined%
\gdef\SetFigFont#1#2#3#4#5{%
  \reset@font\fontsize{#1}{#2pt}%
  \fontfamily{#3}\fontseries{#4}\fontshape{#5}%
  \selectfont}%
\fi\endgroup%
\begin{picture}(2400,2400)(1201,-7561)
\end{picture}%
\\
& \\
$\mathcal{F}_4$ & \hspace{2cm}$\mathcal{F}_4^\perp$
\end{tabular}
\end{center}

\bigskip

\begin{center}
\begin{tabular}{cc}
\begin{picture}(0,0)%
\includegraphics{omega4z1.pstex}%
\end{picture}%
\setlength{\unitlength}{3947sp}%
\begingroup\makeatletter\ifx\SetFigFont\undefined%
\gdef\SetFigFont#1#2#3#4#5{%
  \reset@font\fontsize{#1}{#2pt}%
  \fontfamily{#3}\fontseries{#4}\fontshape{#5}%
  \selectfont}%
\fi\endgroup%
\begin{picture}(2400,2400)(1201,-7561)
\end{picture}%
\hspace{2cm}&\hspace{2cm}\begin{picture}(0,0)%
\includegraphics{orthogonal2.pstex}%
\end{picture}%
\setlength{\unitlength}{3947sp}%
\begingroup\makeatletter\ifx\SetFigFont\undefined%
\gdef\SetFigFont#1#2#3#4#5{%
  \reset@font\fontsize{#1}{#2pt}%
  \fontfamily{#3}\fontseries{#4}\fontshape{#5}%
  \selectfont}%
\fi\endgroup%
\begin{picture}(2400,2400)(1201,-7561)
\end{picture}%
\\
& \\
$\mathcal{F}_4$ & \hspace{2cm}$\mathcal{F}_4^\perp$
\end{tabular}
\end{center}

\bigskip

\begin{center}
\begin{tabular}{cc}
\begin{picture}(0,0)%
\includegraphics{omega4z2.pstex}%
\end{picture}%
\setlength{\unitlength}{3947sp}%
\begingroup\makeatletter\ifx\SetFigFont\undefined%
\gdef\SetFigFont#1#2#3#4#5{%
  \reset@font\fontsize{#1}{#2pt}%
  \fontfamily{#3}\fontseries{#4}\fontshape{#5}%
  \selectfont}%
\fi\endgroup%
\begin{picture}(2400,2400)(1201,-7561)
\end{picture}%
\hspace{2cm}&\hspace{2cm}\begin{picture}(0,0)%
\includegraphics{orthogonal3.pstex}%
\end{picture}%
\setlength{\unitlength}{3947sp}%
\begingroup\makeatletter\ifx\SetFigFont\undefined%
\gdef\SetFigFont#1#2#3#4#5{%
  \reset@font\fontsize{#1}{#2pt}%
  \fontfamily{#3}\fontseries{#4}\fontshape{#5}%
  \selectfont}%
\fi\endgroup%
\begin{picture}(2400,2400)(1201,-7561)
\end{picture}%
\\
& \\
$\mathcal{F}_4$ & \hspace{2cm}$\mathcal{F}_4^\perp$
\end{tabular}
\end{center}

\bigskip

\begin{center}
\begin{tabular}{cc}
\begin{picture}(0,0)%
\includegraphics{omega6.pstex}%
\end{picture}%
\setlength{\unitlength}{3947sp}%
\begingroup\makeatletter\ifx\SetFigFont\undefined%
\gdef\SetFigFont#1#2#3#4#5{%
  \reset@font\fontsize{#1}{#2pt}%
  \fontfamily{#3}\fontseries{#4}\fontshape{#5}%
  \selectfont}%
\fi\endgroup%
\begin{picture}(2400,2400)(1201,-2761)
\end{picture}%
\hspace{2cm}&\hspace{2cm}\begin{picture}(0,0)%
\includegraphics{omega7.pstex}%
\end{picture}%
\setlength{\unitlength}{3947sp}%
\begingroup\makeatletter\ifx\SetFigFont\undefined%
\gdef\SetFigFont#1#2#3#4#5{%
  \reset@font\fontsize{#1}{#2pt}%
  \fontfamily{#3}\fontseries{#4}\fontshape{#5}%
  \selectfont}%
\fi\endgroup%
\begin{picture}(2400,2400)(1201,-7561)
\end{picture}%
\\
& \\
$\mathcal{F}_6$ & \hspace{2cm}$\mathcal{F}_7$
\end{tabular}
\end{center}

\section*{Index des notations}

\bigskip\bigskip

\newlength{\largeur}
\setlength{\largeur}{12.5cm} \addtolength{\largeur}{0cm}
\hspace{-0.5cm}\begin{tabular}{p{1.5cm}p{\largeur}}

\noteA{\mathrm{Sing}(\mathcal{F})}{lieu singulier du feuilletage $\mathcal{F}$}{not0}

\noteA{\mathscr{F}(n;N)}{ensemble des feuilletages de degr\'e $N$ sur l'espace
 projectif $\mathbb{P}^n(\mathbb{C})$}{not7}

\noteA{\mu(\mathcal{F},m)}{nombre de \textsc{Milnor} du feuilletage $\mathcal{F}$
au point singulier $m$}{not1}

\noteA{\langle P,Q\rangle}{id\'eal engendr\'e par les \'el\'ements $P$ et $Q$
de $\mathbb{C}\{x,y\}$}{not1b}

\noteA{\omega_1}{$x^2\mathrm{d}x+y^2(x\mathrm{d}y-y\mathrm{d}x)$}{not2}

\noteA{\omega_2}{$x^2\mathrm{d}x+(x+y^2)(x\mathrm{d}y-y\mathrm{d}x)$}{not3}

\noteA{\omega_3}{$xy\mathrm{d}x+(x^2+y^2)(x\mathrm{d}y-y\mathrm{d}x)$}{not4}

\noteA{\omega_4}{$(x+y^2-x^2y)\mathrm{d}y+x(x+y^2)\mathrm{d}x$}{not5}

\noteA{\omega_5}{$x^2\mathrm{d}y+y^2(x\mathrm{d}y-y\mathrm{d}x)$}{not5b}

\noteA{\mathcal{O}(\mathcal{F})}{orbite 
d'un \'el\'ement $\mathcal{F}$ de $\mathscr{F}(2;2)$ sous l'action 
de~$\mathrm{PGL}_3(\mathbb{C})$}{not8}

\noteA{\mathrm{Iso}(\mathcal{F})}{groupe d'isotropie du 
feuilletage $\mathcal{F}$ d\'efini sur $\mathbb{P}^2(\mathbb{C})$}{not8d}

\noteA{\chi(\mathbb{P}^2(\mathbb{C}))}{alg\`ebre de \textsc{Lie}
des champs de vecteurs holomorphes globaux sur $\mathbb{P}^2(\mathbb{C})$}{not8e}

\noteA{\mathfrak{g}}{alg\`ebre de \textsc{Lie} du groupe alg\'ebrique
$\mathrm{G}$}{not9}

\noteA{\mathcal{R}}{champ radial $x\frac{\partial}{\partial x}+y\frac{\partial}{\partial y}$}{not10}

\noteA{\mathcal{F}_J}{feuilletage d\'ecrit par la $1$-forme $\omega_J$}{not8b}

\noteA{\omega_J}{$(x^2y-1)\mathrm{d}x+(y^2-x^3)\mathrm{d}y$}{not8c}

\noteA{\mathrm{BB}(\mathcal{F};m)}{invariant de \textsc{Baum}-\textsc{Bott} du
feuilletage $\mathcal{F}$ au point $m$}{not11}

\noteA{\mathcal{F}_6}{feuilletage associ\'e au pinceau $x(y+z)+\lambda y(x+z)$}{not14}

\noteA{\mathcal{F}_7}{feuilletage associ\'e au pinceau $xz+\lambda y(y-x)$ }{not15}
\end{tabular}

\vspace*{8mm}

\bibliographystyle{plain}
\bibliography{feuilP2}
\nocite{*}

\end{document}